\documentclass{amsart}
\usepackage{amsmath,amssymb,amsthm,url,xcolor,graphicx}

\newtheorem{thr}{Theorem}
\newtheorem{lem}[thr]{Lemma}
\newtheorem{cor}[thr]{Corollary}

\newtheorem{obs}[thr]{Observation}
\newtheorem{conj}[thr]{Conjecture}
\newtheorem{remr}[thr]{Remark}

\theoremstyle{definition}
\newtheorem{defn}[thr]{Definition}
\newtheorem{conv}[thr]{Convention}
\newtheorem{ex}[thr]{Example}
\newtheorem{quest}[thr]{Question}

\numberwithin{equation}{section}

\def\R{\mathbb{R}}

\def\Rr{{\mathcal R}}

\def\L{{\mathcal L}}

\def\wcc{\operatorname{pc}}
\def\xc{\operatorname{xc}}
\def\dist{\operatorname{dist}}

\def\conv{\operatorname{conv}}
\def\dist{\operatorname{dist}}

\begin{document}

\title{Sublinear extensions of polygons}

\author{Yaroslav Shitov}
\email{yaroslav-shitov@yandex.ru}


\begin{abstract}
Every convex polygon with $n$ vertices is a linear projection of a higher-dimensional polytope with at most $147\,n^{2/3}$ facets.
\end{abstract}

\maketitle

We revisit the following appealing question, which survived a thorough discussion in recent literature without receiving any definitive progress towards its resolution.

\begin{quest}\label{questwcc}
What is the smallest number, $\wcc(n)$, such that every convex polygon with $n$ vertices is a linear projection of some polytope with at most $\wcc(n)$ facets?
\end{quest}

Apart from being a natural question to ask in polyhedral combinatorics, Question~\ref{questwcc} is being studied from the point of view of modern optimization theory and linear algebra. Our approach is purely geometric, but still we need to recall several concepts of optimization theory to explain the relevance of Question~\ref{questwcc}.

\section{Introduction}

Let $P$ be a \textit{polytope}, that is, the convex hull of a finite collection of points in a $d$-dimensional Euclidean space. An \textit{extended formulation} of $P$ is a pair $(Q,\pi)$ consisting of a polytope $Q\subset\R^m$ and a linear projection $\pi:\R^m\to\R^d$ such that $\pi(Q)=P$. The number of facets of such a polytope $Q$ is called the \textit{size} of the formulation $(Q,\pi)$, and the \textit{extension complexity} of  $P$ is the smallest possible size of any extended formulation of $P$. This quantity, further denoted by $\xc(P)$, equals the smallest possible number of inequalities in a representation of any polytope that can be sent to $P$ by a linear projection~\cite{FKPT, Yan}. As we see, Question~\ref{questwcc} is devoted to the worst-case extension complexities of two-dimensional polytopes, so we decided to use the notation `$\wcc$' as an abbreviation of `\textit{polygon complexity}.'

Extended formulations became a prominent technique in optimization after a 1991 paper by Yannakakis~\cite{Yan}, who discussed a possibility of building fast algorithms for hard combinatorial problems by constructing smaller extended formulations to linear programs corresponding to these problems. Several intriguing questions of this nature were solved quite recently, which include the lack of polynomial-size extended formulations for the polytope associated to the traveling salesman problem~\cite{FMPTdW}. Another notable result~\cite{Roth} gives an exponential lower bound for the extension complexity of the matching polytope, which is in contrast to the existence of a polynomial-time algorithm solving the matching problem~\cite{FF}. Other results on the topic include strong lower bounds for the complexities of the polytopes associated to different exact and approximate optimization problems such as max-cut~\cite{CLRS}, independent set~\cite{GJW}, knapsack~\cite{PVV} and many others~\cite{AT}. The extension complexity is also being discussed for polytopes not necessarily arising from particular combinatorial problems, including those with few vertices and facets~\cite{GPPT, GRT, Pad}, the correlation polytope~\cite{KW}, the permutahedron~\cite{Goem}, 0-1 polytopes~\cite{BDP, Roth0}, two-level polytopes~\cite{ACF}, orbitopes~\cite{FK}, Cartesian products and hypersimplices~\cite{GPS} and others. These also include a prominent example of polygons~\cite{B-TN, FRT, ShiCom, VGG}, which this paper is devoted to and which have `prototypical importance' to applications, according to Braun and Pokutta~\cite{BP0}. Also, the study of extended formulations has been developed to reflect the case of semidefinite programming~\cite{FGPRT, GPT, GRT2}, and it may require tools of linear algebra~\cite{Yan}, information theory~\cite{BP}, communication complexity~\cite{FFGT}, quantum learning~\cite{LRS}, tropical mathematics~\cite{mymaprb} and other fields; an interested reader is referred to survey papers~\cite{CCZ, Kai} for further details. 

Yannakakis~\cite{Yan} developed a linear algebraic approach to extended formulations based on \textit{nonnegative matrices}, that is, matrices with nonnegative real entries. The \textit{nonnegative rank} of such a matrix $M$ is the smallest integer $k$ for which $M$ can be written as a sum of $k$ nonnegative rank-one matrices. Now consider a polytope $P$, with $v$ vertices and $f$ facets, defined as a subset of $\R^d$ by the conditions
$$f_\varepsilon(x)=\alpha_\varepsilon,\,\,\,\,g_\varphi(x)\geqslant \beta_\varphi,$$ where the letter $\varepsilon$ indexes a finite set of linear equations, the $\varphi$ runs over the facets of $P$, and the mappings $f_\varepsilon$ and $g_\varphi$ are linear functionals $\R^d\to\R$. A  \textit{slack matrix} of $P$ can be defined by the formula $$S_{\nu\varphi}=g_\varphi(\nu)- \beta_\varphi$$ with rows and columns of $S$ labeled by the vertices and facets of $P$, respectively. It is not hard to see that the rank of $S$ equals $\dim P+1$; Yannakakis~\cite{Yan} proves that the nonnegative rank of $S$ equals the extension complexity of $P$. We note in passing that the nonnegative rank is an important concept of linear algebra in its own right~\cite{Gillis, Vavas} and for theoretical studies in demography, quantum mechanics~\cite{CR}, and statistics~\cite{KRS, MSS}. Also, this notion is relevant for several real-life challenges like signal processing~\cite{GNHS}, text mining~\cite{textmin}, and image processing~\cite{LS}.

\begin{thr}\label{thrYan}\emph{(See~\cite{Yan}.)}
The extension complexity of a polytope $P$ equals the nonnegative rank of any slack matrix of $P$.
\end{thr}


\begin{cor}\label{corlinalg}
For $n\geqslant 3$, the value $\wcc(n)$ equals the largest nonnegative rank of a nonnegative matrix with at most $n$ rows and conventional rank three.
\end{cor}

\begin{proof}
This follows from Theorem~\ref{thrYan} by standard techniques. In particular, we can realize $\wcc(n)$ as the nonnegative rank of an $n\times n$ matrix of rank three by a direct application of this result to a slack matrix of a convex $n$-gon with extension complexity $\wcc(n)$. The fact that the nonnegative rank of any $n\times m$ nonnegative matrix $A$ with rank three cannot exceed $\wcc(n)$ also follows from Theorem~\ref{thrYan}, as outlined in the proof of Theorem~3.1 in~\cite{ShiCom}. To give a sketch of the proof, we consider such a matrix $A$, and we denote by $\Delta\subset\R^n$ the simplex consisting of nonnegative vectors whose coordinates sum to $1$. Since $\Delta$ has $n$ facets, the intersection of $\Delta$ with the column space of $A$ is a polygon $P$ with $k\leqslant n$ vertices. Denoting by $S$ the matrix of the column coordinate vectors of the vertices of $P$, we have $A=SB$ with $B$ nonnegative. Since $S$ is a slack matrix of $P$, it has nonnegative rank at most $\wcc(n)$, and so does $A$.
\end{proof}



\section{Polygons on a plane}\label{secresold}

Question~\ref{questwcc} takes an important place in the study of extended formulations. A detailed account of this problem was carried out by the authors of~\cite{FRT}, who adopted the dimension counting method and gave an $\Omega(\sqrt{n})$ lower bound for the extension complexity of a \textit{generic} convex $n$-gon. More presicely, they obtained the inequalities
\begin{equation}\label{eq1}
\sqrt{2n}\leqslant\wcc(n)\leqslant n,
\end{equation}
in which the upper bound is trivial. Question~\ref{questwcc} and its analogues were studied in at least a dozen of papers~\cite{BL, GG, GPT, GRT2, Hr, LC, Pad, PP, PW, ShiCom, mymaprb, my11, my12, VGGT} and mentioned in several dozens of other works, including a number of textbooks, surveys and highly cited research papers~\cite{BP, CCZ, FKPT, Moitra3, RY, Roth0, Venkat}. This question was discussed at open problem sessions of computer science conferences~\cite{Dag13, Dag15} and online in popular media~\cite{Garden, MO}, but it has seen no progress on either bound except for improving the constant factors in front of $\sqrt{n}$ and $n$. Concerning the asymptotic behavior of $\wcc(n)$, the authors of the initial paper~\cite{FRT} and most other experts seemed to expect that the upper bound in~\eqref{eq1} is closer to the actual value. 
Braun and Pokutta~\cite{BP} point out the importance of the following question and mention it alongside a list of major problems in the theory of extended formulations. 

\begin{quest}\label{questwccn}
Do we have $\wcc(n)\geqslant\varepsilon n$ for some fixed $\varepsilon>0$? 
\end{quest}

This question appeared in online media such as the \textit{Open Problem Garden}~\cite{Garden} and \textit{MathOverflow}~\cite{MO}. An affirmative answer to Question~\ref{questwccn} was conjectured in~\cite{BL} and claimed to be proven in~\cite{LC}. However, Hrube\v{s}~\cite{Hr} showed that the argument of~\cite{LC} is flawed and formulated a problem equivalent to Question~\ref{questwccn}; a further version of this question but specified to $\varepsilon=0.5$ appears in~\cite{VGGT}. The authors of~\cite{FRT} expected an affirmative answer to another version obtained by replacing $\varepsilon n$ with $\tilde{\Omega}(n)$ in the formulation of Question~\ref{questwccn}, which would mean that the upper bound in~\eqref{eq1} is tight up to logarithmic factors. Contrary to these expectations, we show that the lower bound in~\eqref{eq1} is closer to the actual value of $\wcc(n)$ in the asymptotical sense. 

\begin{thr}\label{thrmain}
We have $\wcc(n)\leqslant 147\,n^{2/3}$ for all $n\geqslant 3$.
\end{thr}

We close the preliminary part of our paper with a more detailed description of the progress achieved on Questions~\ref{questwcc} and~\ref{questwccn} before the work we present. As far as we know, Beasley and Laffey~\cite{BL} were the first to consider this question; they followed the linear algebraic approach in terms of Corollary~\ref{corlinalg}. They proved that a special family of nonnegative $n\times n$ rank-three matrices, called \textit{Euclidean distance matrices}, have nonnegative rank at least $\log n$, and conjectured that the maximal value of this rank is $n$. In our notation, their result stated that $\wcc(n)\geqslant\log n$, and the conjecture was $\wcc(n)=n$. A subsequent paper~\cite{LC} claimed to prove this conjecture, but later works~\cite{GG, Hr} pointed out a flaw in their proof. The authors of~\cite{GG} saved a part of the argument in~\cite{LC} and showed that the so-called \textit{restricted nonnegative rank} of Euclidean distance matrices equals $n$. They reiterated the question as to whether the equality $\wcc(n)=n$ holds for general $n$ and proved it for $n$ not exceeding five. The $n=6$ case of this question was treated in the affirmative in~\cite{GPT}, but already for $n=7$ the answer turned out to be negative --- the equality $\wcc(7)=6$ was proved in~\cite{PP, ShiCom}. As a corollary of the latter result, the authors of~\cite{PP, ShiCom} get the estimate $\wcc(n)\leqslant\lceil 6n/7\rceil$ for all $n$, which stood as the best known upper bound on $\wcc(n)$ to this date. A more thorough description of the $n=7$ case lead the authors of~\cite{PW} to a family of polytopes whose optimal extended formulations require the use of \textit{hidden vertices}. The tropical approach to Questions~\ref{questwcc} and~\ref{questwccn} does not look promising; the author of~\cite{my11, my12} gave new combinatorial proofs that $\wcc(6)=6$ and $\wcc(8)=7$, but both of these equalities were known at the time because an improved counting argument of Padrol~\cite{Pad} confirmed the inequality $\wcc(n)\geqslant2\sqrt{2n-2}-1$ for all $n$. An algorithmic treatment of these questions gave the authors of~\cite{VGGT} some important information on polygons whose vertices are randomly chosen on a circle, and they went on to conjecture that $\wcc(n)\leqslant 0.5 n+3$ for all $n$. The paper~\cite{GRT2} showed that any convex $n$-gon admits a \textit{semidefinite} extended formulation of size at most $4\lceil n/6\rceil$, but this did not imply any progress on our $\wcc$ function because the semidefinite approach has more expressive power in comparison to the extended formulations as in the setting of this paper. Nevertheless, Question~\ref{questwccn} remained open even in the case of stronger, semidefinite, extended formulations.

\section{Preliminaries}

This paper is devoted to the proof of Theorem~\ref{thrmain}. In this section, we collect several standard notational conventions and basic results used in the course of our discussion, and we give a short overview of our approach.

We work in the Euclidean real space $\R^d$ with the conventional Euclidean metric. The notation $\dist(u,v)$ stands for the distance between two points $u, v\in\R^d$, and the distance $\dist(U,V)$ between two non-empty subsets $U, V\subset\R^d$ is the infimum of all possible distances $\dist(u, v)$ with $u\in U$ and $v\in V$. The \textit{interior} of a set $V\subset\R^d$ is understood in terms of the topology inherited from the Euclidean metric, and the \textit{relative interior} of a segment $\sigma\subset\R^d$ is the same segment but with endpoints removed. A set $V\subset\R^d$ is called \textit{convex} if for any pair of distinct points $u, v\in V$, the segment between $u$ and $v$ is contained in $V$. The convex hull of a set $S\subset\R^d$ is denoted by $\conv S$, and it is defined as the intersection of all convex sets in $\R^d$ containing $S$. We define a polytope in $\R^d$ and a polygon in $\R^2$ as the convex hull of a finite collection of points, that is, we do not work with non-convex generalizations of polytopes here. If $u, v\in\R^d$ are distinct points, then we write $\overrightarrow{uv}$ for the \textit{oriented segment} with the beginning at $u$ and ending at $v$. We note that this is different from saying that $\overrightarrow{uv}$ is a \textit{vector}, which would rather mean the class of all segments with the same length and direction.

A significant part of our considerations goes over the two-dimensional space, so we give several definitions that we mostly use in this setting. A \textit{ray} $\rho\subset \R^2$ is a closed convex unbounded subset of a straight line; if $a$ is the apex and $b$ is a non-apex point on $\rho$, then we say that $\rho$ goes \textit{from} $a$ \textit{towards} $b$. If $r, s$ are two oriented segments or rays, then we write $\angle(r, s)$ for the measure of the angle between them, and we also write $\angle uvw$ for $\angle(\overrightarrow{vu}, \overrightarrow{vw})$ whenever $u,v,w$ are points.

\begin{obs}\label{obsaddang}
If 
$a_1, a_2, a_3$ are oriented segments or rays, then $$\angle(a_1, a_2)+\angle(a_2, a_3)\geqslant \angle(a_1, a_3).$$
\end{obs}

If $P\subset\R^2$ is a polygon, then the \textit{turning angle} at a vertex $v$ is $\pi-\angle v_-vv_+$, where $v_-$ and $v_+$ are two vertices adjacent to $v$ in $P$. The \textit{turning angle} of an edge $e$ of a polygon is the sum of the turning angles at the two endpoints of $e$.

\begin{defn}\label{defnwedge}
We write $u\wedge v$ to denote

\noindent (1) the straight line connecting $u, v$ if they are distinct points,

\noindent (2) the intersection of $u, v$ if they are non-collinear straight lines on a plane.
\end{defn}

We proceed with two basic results on the extension complexity.

\begin{lem}\label{lemeasy2}
Let $P$ and $Q$ be polytopes in $\R^d$ each of which is different from a single point. Then $\xc(\operatorname{conv} P\cup Q)\leqslant\xc(P)+\xc(Q)$.
\end{lem}

\begin{proof}
See Proposition~3.1.1 in~\cite{WeltgePhD} and a similar Proposition~2.8(6) in~\cite{GPT}.
\end{proof}

\begin{obs}\label{lemeasy1}
Let $P\subset\R^d$ be a polytope, and let $H\subset\R^d$ be a closed half-space. If the intersection $P\cap H$ is non-empty, then $\xc(P\cap H)\leqslant \xc(P)+1$.
\end{obs}

\begin{proof}
Let $Q\subset\R^d\times\R^\alpha$ be a polytope with $\xc(P)$ facets such that $\pi(Q)=P$, where $\pi$ is the projection onto the first $d$ coordinates. Then the polytope $Q'$ defined as $Q\cap\left(H\times\R^\alpha\right)$ has at most $\xc(P)+1$ facets, and we have $\pi(Q')=P\cap H$.
\end{proof}

The general strategy of our approach is as follows. The forthcoming Sections~\ref{secacute}--\ref{secreduce} reduce the question of constructing a small extended formulation of a polygon to a specific problem in plane geometry, which comes from a certain special class of three-dimensional extensions. This class is introduced in Section~\ref{secacute} under the name of \textit{acute polyhedra}, and their Schlegel diagrams are investigated in Section~\ref{secdiag}. Section~\ref{secreduce} contains Theorem~\ref{lemacutemult}, which is the main result of the first part of the paper, and it explains how to glue several extensions coming from acute polytopes together to get a higher dimensional extended formulation of an initial polygon. 

In view of Lemma~\ref{lemeasy2}, the progress on Question~\ref{questwcc} can be made by extracting a sufficiently large subset $S$ of the vertices of a given polygon for which the extension complexity of $\conv S$ is small, so the subsequent Sections~\ref{secthin}--\ref{secslant} are devoted to extracting such a subset which would be well behaved with respect to the methods developed in Sections~\ref{secacute}--\ref{secreduce}. More precisely, Sections~\ref{secthin} and~\ref{secgoodseq} give further notational conventions needed to get an appropriate two-dimensional description of polygons that allow the use of Theorem~\ref{lemacutemult}. In Section~\ref{secextr}, we introduce a further auxiliary notion of when the sequence of the vertices of a given polygon is \textit{slanted}, and we explain that a polygon without large slanted subsequences should have relatively small extension complexity. The problem reduces to proving good upper bounds of polygons with slanted vertex sequences, which is done in the technical Sections~\ref{secdecr} and~\ref{secslant}. Our proof is completed in Section~\ref{seccompl}, and we give several additional remarks and suggestions on further research on this topic in Section~\ref{secrem}.

\section{Acute polyhedra}\label{secacute}

In this section, we introduce the class of \textit{acute} polytopes, which are useful for constructing small extended formulations. We recall that, for three-dimensional polyhedra, the angle between two faces $A$ and $B$ with a common edge $e$ equals the angle between two oriented segments that lie on $A$ and $B$, respectively, have their beginnings on $e$, and are orthogonal to $e$.

\begin{defn}\label{defacute}
Assume $P\subset\R^d$ is a polytope and $B$ is one of its facets. If any other facet $F\neq B$ of $P$ satisfies the conditions that

(a) $\dim B\cap F=\dim P-2$ and

(b) the angle between $B$ and $F$ is acute,

\noindent then $P$ is called an \textit{acute} polytope with \textit{base} facet $B$.
\end{defn}

\begin{remr}\label{remrdefacute}
Our notation is suggested by the fact that an acute triangle satisfies the above conditions whichever side of it we choose as the base. Right and obtuse triangles can also become 'acute polytopes' in our sense, but the only possible bases are their longest sides. For $n>3$, no $n$-gon on the plane can satisfy Definition~\ref{defacute} because of the condition~(a). In our further considerations, we restrict our attention to three-dimensional acute polytopes, and we call them \textit{acute polyhedra}.
\end{remr}

One example of an acute polyhedron is given on Figure~\ref{figacute} below.

\begin{obs}\label{lemacute1}
Let $P\subset\R^3$ be an acute polyhedron with base facet $B$. Assume that $B$ lies on the plane $\{z=0\}$ and its edge-defining inequalities are $\{a_i x + b_i y+c_i\geqslant 0\}$ with $i=1, \ldots, k$. Then the non-base facets of $P$ are defined by 
$$a_i x + b_i y+c_i\geqslant\varepsilon_i z,$$
where $\varepsilon_1,\ldots,\varepsilon_k$ are non-zero numbers with the same sign.
\end{obs}

\begin{proof}
According to the part (a) of Definition~\ref{defacute}, every non-base facet of $P$ passes through some edge of $B$. Any such facet should have the inequality as in the displayed formula, and the part (b) of Definition~\ref{defacute} shows that the $\varepsilon$'s should be taken positive if $P\subset\{z\geqslant 0\}$ and negative if $P\subset\{z\leqslant 0\}$.
\end{proof}

We need one more concept to be used in further considerations.

\begin{defn}\label{defmainedge}
We say that an edge $e$ of an acute polyhedron with base $B$ is \textit{main} if exactly one of the endpoints of $e$ lies on $B$.
\end{defn}

\begin{lem}\label{lemunique}
Any base vertex of an acute polyhedron belongs to a unique main edge.
\end{lem}

\begin{proof}
It is clear from Definition~\ref{defmainedge} that a base vertex $v$ is adjacent to two  other base vertices, so if the current statement was false, there would be at least four vertices adjacent to $v$. This means that there would be at least four faces containing $v$, in which case we can find three such faces $\varphi_1,\varphi_2,\varphi_3$ none of which is the base. The part (a) of Definition~\ref{defacute} shows that these faces pass through base edges $e_1, e_2, e_3$, respectively. If any such $e_i$ did not contain $v$, the face $\varphi_i$ would pass through three non-collinear base vertices, which is a contradiction because $\varphi_i$ is not the base. Therefore, the base edges $e_1, e_2, e_3$ pass through $v$, which means that two of them coincide. This implies that the base and two of the faces $\varphi_1, \varphi_2, \varphi_3$ have a common edge, which is impossible for three-dimensional polyhedra. 
\end{proof}

\begin{lem}\label{lemdir}
Let $V\subset\R^2$ be a convex polygon with vertices $v_0, v_1, \ldots, v_n$, and let $y_1, \ldots, y_n$ be a family of inner points of $V$. Then there is an acute polyhedron $P\subset\R^3$ with base $V$ such that, for any $i\in\{1,\ldots,n\}$, the image of the main edge passing from $v_i$ under the orthogonal projection of $P$ onto $V$ is collinear to $v_i\wedge y_i$.
\end{lem}

\begin{proof}
We embed the plane containing $V$ into $\R^3$ and choose one of the two resulting half-spaces $U$ to be called \textit{upper}.  For all $i\in\{1,\ldots,n\}$, we define $H_i$ as the plane that passes through $v_i, y_i$ and is orthogonal to the plane of $V$. Also, we denote the edges of $V$ consecutively as $e_0,\ldots,e_n$, and we assume that $e_j$ has $v_{j}$ and $v_{j+1}$ as endpoints, where $v_{n+1}$ stands for $v_0$ and $j$ ranges in $0, \ldots, n$. 

In the following construction of an acute polyhedron $P\subset U$, the notation $F_j$ stands for the upper half-plane that contains $e_j$ and the corresponding non-base facet of $P$. We choose $F_0$ as an arbitrary upper half-plane passing through $e_0$ and having an acute angle with $V$, and, by the formulation of the lemma, we conclude that the main edge passing through $v_1$ should be a part of the ray $F_0\cap H_1$, which allows us to define $F_1$ as the upper half-plane containing $e_1$ and $F_0\cap H_1$. We proceed by the induction, and we define $F_i$ as the upper half-plane passing through $e_i$ and $F_{i-1}\cap H_i$ for all $i\in\{1,\ldots,n\}$. It remains to note that $F_i$ forms an acute angle with $V$ because it contains a non-apex point on $F_{i-1}\cap H_i$ which appears inside $V$ under the orthogonal projection onto the plane containing $V$.
\end{proof}

\section{Acute diagrams}\label{secdiag}

Our further considerations need a more accurate characterization of acute polyhedra in terms of their  \textit{Schlegel diagrams}~\cite{GrunBook}. In this section, we describe the class of planar straight-line graphs that can arise as the diagrams of acute polyhedra.

\begin{obs}\label{obsproj}
Let $P\subset\R^3$ be an acute polyhedron; consider the orthogonal projection $\pi$ onto the plane containing its base $B$. Then (a) the non-base points of $P$ project into the interior of the base, and (b) the mapping $\pi$ is injective on the non-base boundary points of $P$.
\end{obs}

\begin{proof}
The conclusion (a) is immediate from the part (b) of Definition~\ref{defacute}. In order to prove (b), we consider two non-base points $u\neq v$ on the boundary of $P$ with $\pi(u)=\pi(v)$. Since $B$ is a facet, one of the points $u, v$ should lie between the plane of $B$ and the other point. If this middle point is $u$, then by the item (a) it lies in the interior of $\conv B\cup\{v\}$ and hence in the interior of $P$, which is a contradiction.
\end{proof}

\begin{remr}\label{remrdiag}
In other words, the projection $\pi$ as in Observation~\ref{obsproj} gives a \textit{Schlegel diagram} of $P$. This diagram is to be called the \textit{acute diagram} of $P$ \textit{relative to} the base facet $B$, or just the acute diagram of $P$ if the choice of the base is clear from context. We also say that $P$ is an \textit{acute lifting} of the corresponding acute diagram.
\end{remr}

\begin{figure}[h]
\includegraphics[width=10cm]{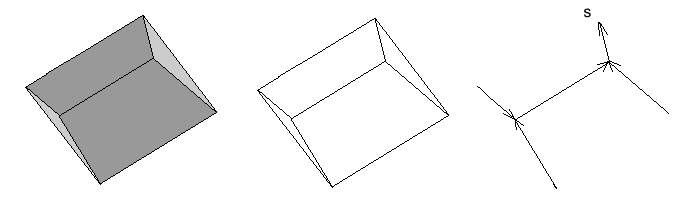}
\caption{An acute polytope, its acute diagram, and an $s$-orientation.}
\label{figacute}
\end{figure}

Figure~\ref{figacute} gives an example of an acute polytope and the corresponding diagram.

\begin{defn}\label{defnpslg}
A \textit{planar straight-line graph} $\Delta$ is a topological graph in $\R^2$ whose \textit{arcs} are non-crossing line segments, that is, any pair of arcs are either disjoint or have a common \textit{node}. The convex hull of $\Delta$ is to be called the \textit{base} of $\Delta$, and the vertices of the base are called the \textit{base nodes} of $\Delta$.
\end{defn}

We say that line segments $a, b, c\subset\mathbb{R}^2$ are \textit{concurrent} if the straight lines containing $a$, $b$, $c$ are concurrent, that is, are either parallel or have a common point.

\begin{lem}\label{lemproj}
Let $\Delta$ be the acute diagram of an acute polyhedron $P$ with base $B$. Then $\Delta$ is a planar straight-line graph such that

(o) the base of $\Delta$ is $B$,

(i) every node of $\Delta$ has degree at least three,

(ii) the non-base nodes of $\Delta$ lie in the interior of the base,

(iii) every edge of the base is an arc of $\Delta$,

(iv) every bounded face $f$ of $\Delta$ contains exactly one arc $e_f$ of the base,

(v) if a non-base arc $e$ of $\Delta$ separates faces $f, g$, then $e, e_f, e_g$ are concurrent.
\end{lem}

\begin{proof}
The planarity of $\Delta$ and the item (i) follow from the part (b) of Observation~\ref{obsproj}, and the items (o, ii, iii) can be obtained from the part (a) of Observation~\ref{obsproj}. To prove the item (iv), we note that any face of $P$ contains at least one edge of the base by the part (a) of Definition~\ref{defacute}, and a face with at least two such edges should be the base itself.

Now we proceed with the item (v). Let $\varphi,\psi$ be two faces of $P$ with a common edge $\varepsilon$. We denote by $\overline{B}$ the base plane of $P$, and by $\Phi, \Psi$ the two planes containing $\varphi,\psi$, respectively. Then $\Phi\cap\Psi$ is the straight line containing $\varepsilon$, and its intersection with the base plane is $\Phi\cap\Psi\cap \overline{B}$, which is the same as the intersection of the two straight lines $\Phi\cap \overline{B}$, $\Psi\cap \overline{B}$ containing the base edges of $\varphi$, $\psi$.
\end{proof}

\begin{remr}\label{remrorient}
Let $\Delta$ be a planar straight-line graph satisfying the conditions (i)--(v) as in Lemma~\ref{lemproj}, and let $s$ be an arbitrary base node of $\Delta$. We define $\Delta^{\operatorname{in}}$ as the graph obtained by removing all the base arcs from $\Delta$. According to the item (iv) of Lemma~\ref{lemproj}, this graph is a tree, so we can define the $s$\textit{-orientation} by saying that an arc $\conv\{u,v\}$ of $\Delta^{\operatorname{in}}$ is oriented from $u$ to $v$ if and only if $\rho(u,s)>\rho(v,s)$, where $\rho$ is the graph-theoretical geodesic distance function on $\Delta^{\operatorname{in}}$.
\end{remr}

Now we need to explain how to lift an acute diagram into an acute polytope.

\begin{obs}\label{obsoutdeg}
In the $s$-orientation as in Remark~\ref{remrorient}, every node different from $s$ has outdegree $1$. The node $s$ is a sink, that is, it has outdegree $0$.
\end{obs}

\begin{proof}
Since $\Delta^{\operatorname{in}}$ is a tree, the shortest path between any pair of nodes is unique, which proves the first statement. The conclusion on the outdegree of $s$ is trivial.
\end{proof}

\begin{lem}\label{lemextra}
Let $u$ be an inner node of the $s$-orientation of $\Delta$ as in Remark~\ref{remrorient}. Let $\Lambda(u)$ be the set of all base nodes from which $u$ can be reached by oriented paths. Then $\Lambda(u)$ is a set of consecutive nodes of the base.
\end{lem}

\begin{proof}
For $v_1, v_2$ in $\Lambda(u)$, we can construct oriented paths $\pi_1, \pi_2$ beginning at $v_1, v_2$, respectively, and ending at $u$. The union of $\pi_1$ and $\pi_2$ splits $\Delta$ into two connected components, and for any base node $x$ outside the component of the sink, any oriented path from $x$ to $s$ passes through $u$.
\end{proof}

\begin{lem}\label{lemmainint0}
Any planar straight-line graph $\Delta$ satisfying the conditions (i)--(v) as in Lemma~\ref{lemproj} contains a triangle formed by two main arcs and one base arc $e$. This arc $e$ can be chosen so that it does not contain a base node $s$ fixed in advance.
\end{lem}

\begin{proof}
We choose a node $x$ furthest from $s$ in terms of the graph-theoretic geodesic distance on $\Delta^{\operatorname{in}}$. Since $\Delta^{\operatorname{in}}$ is a tree, this $x$ should be its leaf, that is, a base node with a unique main arc $\alpha$ leading from $x$ to some inner node $y$. By the condition (i) in Lemma~\ref{lemproj} and Observation~\ref{obsoutdeg}, we see that $y$ is adjacent to at least two nodes $x, x'$ at the maximal distance from $s$. According to Lemma~\ref{lemextra}, the base nodes from which $y$ can be reached are consecutive, and any such node should be connected with $y$ by an arc because of the distance maximality assumption.
\end{proof}

\begin{lem}\label{lemmainint}
Let $\Delta$ be a planar straight-line graph satisfying the conditions (i)--(v) as in Lemma~\ref{lemproj}. If the base of $\Delta$ is not a triangle, then there is a base arc $e$ with turning angle at most $\pi$ which forms a triangle together with two main arcs adjacent to the endpoints of $e$. If, additionally, the base of $\Delta$ is not a trapezoid, then this arc $e$ can be chosen to have the turning angle strictly less than $\pi$.
\end{lem}

\begin{proof}
Let us call \textit{bad} those arcs of the base of $\Delta$ that have turning angle $\pi$ or greater. If all bad arcs are adjacent to a common node $s$, then the proof is completed by the application of Lemma~\ref{lemmainint0}. If there are two non-adjacent bad arcs, then the total of the turning angles at all four of their endpoints is $2\pi$, so the base of $\Delta$ has no other nodes except these four, so its further analysis is straightforward.
\end{proof}

\begin{lem}\label{lemlift}
A planar straight-line graph $\Delta$ satisfying the conditions (i)--(v) as in Lemma~\ref{lemproj} can be obtained as an acute diagram of some acute polyhedron $P$.
\end{lem}

\begin{proof}
We work by the induction on $k$, the number of nodes of the base of $\Delta$. The base of the induction comes from the two cases in which Lemma~\ref{lemmainint} is not applicable, namely, the triangle and trapezoid bases. These cases are easy because the triangle bases come from triangular pyramids, and the trapezoid case comes from cutting the triangular prism from the two sides as in Figure~\ref{figacute}.

We proceed with an application of Lemma~\ref{lemmainint}, and we find a base arc $e=\conv\{x_1,x_2\}$ with turning angle less than $\pi$ such that both $x_1,x_2$ are adjacent to the same inner node $y$. Let $\varepsilon_1, \varepsilon_2$ be the base arcs that are different from $e$ and contain $x_1, x_2$, respectively. Now to construct a new graph $\Delta'$ from $\Delta$, we

(1) continue $\varepsilon_1, \varepsilon_2$ to the intersection point $x$,

(2) connect $x$ to $y$, and erase the arcs $\conv\{x_1, y\}$ and $\conv\{x_2, y\}$,

\noindent and we do not call $x_1$ and $x_2$ the nodes of $\Delta'$. If the degree of $y$ became less than three after the transformations (1)--(2), then we do not call it a node as well.

\begin{figure}[h]
\includegraphics[width=10cm]{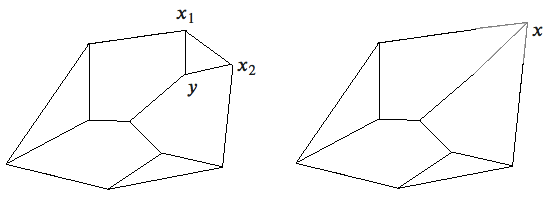}
\caption{The inductive step in the proof of Lemma~\ref{lemlift}.}
\label{figacute2}
\end{figure}

It is straightforward to see that the graph $\Delta'$ is planar and straight-line, and it satisfies the conditions (i)--(v) as in Lemma~\ref{lemproj}. Using the inductive assumption, we build an acute lifting $P'$ of $\Delta'$, and it remains to cut off a piece of $P'$ by the plane passing through $x_1, x_2, y_0$, where $y_0$ is the point on the main edge of $P'$ passing from $x$ that maps to $y$ under the orthogonal projection of $P'$ onto its base.
\end{proof}


We finalize the section with a more delicate general property of acute diagrams. 

\begin{cor}\label{coracute}
Let $V$ be a convex polygon with vertices $v_0, v_1, \ldots, v_n$, and let $y_1, \ldots, y_n$ be a family of inner points of $V$. Then there is an acute diagram with base $V$ such that, for any $i\in\{1,\ldots,n\}$, the main edge from $v_i$ lies on $v_i\wedge y_i$.
\end{cor}

\begin{proof}
This a reformulation of Lemma~\ref{lemdir}.
\end{proof}

\section{Glueing acute extensions together}\label{secreduce}



Now we are going to employ the approach of acute polyhedra in a construction of small extended formulations. We recall that a \textit{polyhedral cone} in $\R^d$ is the convex hull of a finite collection of rays passing from the same apex $O$. Such a cone is called \textit{pointed} if it does not contain a straight line.



\begin{obs}\label{obscut}
Let $C\subset\R^q$ be a pointed cone defined as the convex hull of its extreme rays $l_1,\ldots,l_k$ passing from the apex $O$. Let $H\subset\R^q$ be a closed half-space containing an unbounded part of every $l_1,\ldots,l_k$ such that the boundary of $H$ intersects $l_i$ at a unique point $v_i$ for all $i\in\{1,\ldots,k\}$. Also, let $\pi:\R^q\to\R^d$ be a linear projection. Then $\pi(H\cap C)$ is the convex hull of the union of the rays passing from $\pi(v_i)$ towards the direction of $\pi(l_i)$ taken over all $i=1,\ldots,k$.
\end{obs}

\begin{proof}
The ray $H\cap l_i$ passes from $v_i$ in the direction of $l_i$, and the set $H\cap C$ is the convex hull of the union of all $H\cap l_i$ over $i=1,\ldots,k$. Since the mappings $\pi$ and $\conv$ on the subsets of $\R^q$ commute, we get the desired result.
\end{proof}

Now we are ready to explain how to glue several three-dimensional acute extensions together to get quite a small higher-dimensional extended formulation.

\begin{figure}[h]
\includegraphics[width=12cm]{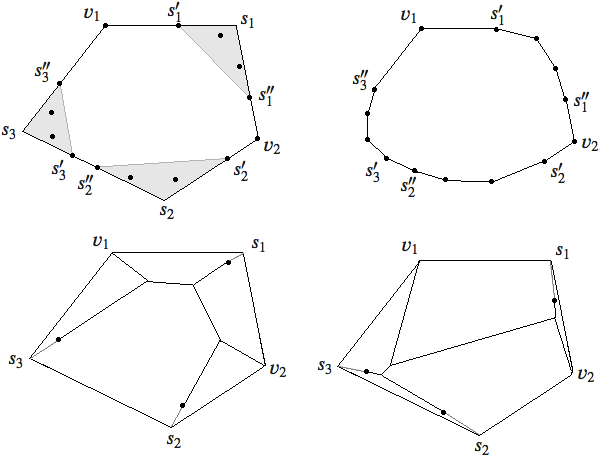}
\caption{An application of Theorem~\ref{lemacutemult}: Two acute diagrams confirm that the $14$-gon on the right has $\xc$ at most $10$.}
\label{figacuteglue}
\end{figure}

\begin{thr}\label{lemacutemult}
Let $P$ be a convex polygon on a plane with a distinguished non-empty set $S$ of the set $V$ of all vertices. Assume that, for any $s\in S$, one chooses two vertices $s', s''$ on two edges of $P$ adjacent to $s$, and one picks a set $\{s^1,\ldots,s^\delta\}$ consisting of $\delta\geqslant 1$ new points in the interior of the triangle $T_s=\operatorname{conv}\{s,s',s''\}$. Assume that the choice of $s', s''$ was made so that the triangles $T_s$ are disjoint for different $s$. Assume that, for any $i\in\{1,\ldots,\delta\}$, there exists an acute diagram $D^i$ with base face $P$ such that, for any $s\in S$, the segment between $s$ and $s^i$ is a subset of the main edge of $D^i$ passing from $s$. Then the convex hull of 
\begin{equation}\label{eqbigconvexhull}
\left(\bigcup\limits_{s\in S}\{s', s'', s^1,\ldots,s^\delta\}\right)\cup\left(V\setminus S\right)
\end{equation}
is a polygon of extension complexity not exceeding $|V|+|S|+\delta$.
\end{thr}

\begin{proof}
We are going to construct a polytope $\mathcal{P}'$ with at most $|V|+|S|+\delta$ facets in $\R^2\oplus\R^\delta=\{(x,y,z^1,\ldots,z^\delta)\}$ such that the projection of $\mathcal{P}'$ onto the first two coordinates is the desired set, that is, the convex hull of~\eqref{eqbigconvexhull}. Let $a_e x+b_e y+c_e\geqslant 0$ be a defining inequality corresponding to an edge $e$ of the polygon $P$. According to Observation~\ref{lemacute1}, the lifting of the diagram $D^i$ gives an acute polytope $A^i\subset\{(x,y,z^i)\}$ defined by the conditions $z^i\geqslant0$ and $$a_e x + b_e y+c_e\geqslant\varepsilon^i_{e} z^i$$ with the $\varepsilon$'s being positive real numbers. We define $\mathcal{P}\subset\R^2\oplus\R^\delta$ by the conditions
\begin{equation}\label{equation1234z}
z^i\geqslant 0,
\end{equation}
\begin{equation}
\label{equation1234e}
a_e x + b_e y+c_e\geqslant\varepsilon^1_{e} z^1+\ldots+\varepsilon^\delta_{e} z^\delta,
\end{equation}
where $i$ ranges in $\{1,\ldots,\delta\}$, and $e$ runs over all edges of $P$. It is clear that $\dim\mathcal{P}=\delta+2$, and we use the part (a) of Observation~\ref{obsproj} to see that $\pi(\mathcal{P})=P$, where $\pi:\R^2\oplus\R^d\to\R^2$ is the projection onto $\{(x,y)\}$. Now we pick any $s=(x_s, y_s)\in S$, and we see that the point $\sigma=(x_s, y_s, 0,\ldots,0)$ is the vertex of $\mathcal{P}$, which fact is easy because $\sigma$ is the vertex of the face of $\mathcal{P}$ defined by $z^1=\ldots=z^\delta=0$. Since there are exactly $\delta+2$ facets of $\mathcal{P}$ intersecting at $\sigma$, namely all those in~\eqref{equation1234z} and those in~\eqref{equation1234e} corresponding to $s\in e$, we conclude that the vertex figure of $\sigma$ is a simplex. In particular, the cone corresponding to $\sigma$ is defined by $\delta+2$ extreme rays, two of which correspond to the two edges of $P$ adjacent to $s$ and the other $\delta$ are~\footnote{We note that a ray of the form~\eqref{equationdir1} contains an edge of $\mathcal{P}$ because it contains, by the definition of an acute diagram, an edge of the corresponding acute polyhedron $A^i$, and this polyhedron is a face of $\mathcal{P}$ defined by the conjunction of the equations $z^j=0$ with $j\neq i$.}
\begin{equation}\label{equationdir1}
(x_s+\alpha_s^i t,y_s+\beta_s^it,0,\ldots,0,t,0,\ldots,0),\,\,\,\,t\geqslant0,
\end{equation}
where we have $i-1$ zeros in front of $t$, and $(\alpha_s^i,\beta_s^i)$ is a vector pointing to the direction from $s$ towards $s^i$. Now we define the half-space $H_s\subset\R^2\oplus\R^\delta$ which

(1) does not contain $\sigma$,

(2) contains $(x_{s'}, y_{s'}, 0,\ldots,0)$ and $(x_{s''}, y_{s''}, 0,\ldots,0)$,

(3) for any $i=1,\ldots,\delta$, it contains the point on the ray~\eqref{equationdir1} that projects to $s^i$.

\noindent We can define $\mathcal{P}'$ as the intersection of $\mathcal{P}$ and all $H_s$ with $s\in S$, and we conclude that $\pi(\mathcal{P})$ is the desired polygon by applying Observation~\ref{obscut}.
\end{proof}

\section{Thin sequences}\label{secthin}

We proceed with several further notational conventions. A sequence $v=(v_1,\ldots,v_n)$ of distinct points on a plane is called \textit{correct} if these points are the vertices of their convex hull $P$ and the segment between any pair of consecutive points in $v$ is an edge of a polygon $P$. A correct sequence is called a \textit{cw-sequence} if the order of the vertices is clockwise, and it is called a \textit{ccw-sequence} otherwise.

\begin{defn}\label{defthin}
Let $\alpha\in(0,\pi)$ and $n\geqslant 3$. A sequence $v=(v_1,\ldots,v_n)$ is called $\alpha$\textit{-thin} if $v$ is correct and, additionally, the turning angle of the edge $\conv\{v_1, v_n\}$ in the polygon $\conv v$ is greater than $2\pi-\alpha$, that is, $$\angle v_n v_1 v_2+\angle v_1v_nv_{n-1}<\alpha.$$ We say that $v$ is \textit{thin} if it is $\alpha$-thin for some $\alpha\in(0,\pi)$.
\end{defn}

In the following easy observation, we use the total turning angle as a parameter of the boundary of a given convex polygon $P$ to split the vertices of $P$ into the union of several thin sequences.

\begin{figure}[h]
\includegraphics[width=8cm]{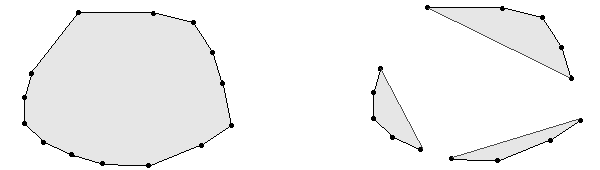}
\caption{Splitting the vertices of a polygon into thin sequences.}
\label{figthinsplit}
\end{figure}

\begin{obs}\label{lemextrthin}
Let $P$ be a convex polygon on $n$ vertices, and let $q\geqslant 3$ be an integer. Then the vertices of $P$ can be partitioned into at most $q$ sets each of which is either a point or a pair of points, or else a set that forms a $2\pi/q$-thin sequence.
\end{obs}

\begin{proof}
We consider a consecutive enumeration $e_1,\ldots,e_n$ of the relative interiors of the edges of $P$. We have $$0=\tau(e_1)<\tau(e_2)<\ldots<\tau(e_n)<2\pi,$$ where $\tau(e_i)$ is the sum of the turning angles at all the vertices that lie between $e_1$ and $e_i$ in the sense of our enumeration. We also denote by $e_{n+1}$ a copy of $e_1$, and we assume $\tau(e_{n+1})=2\pi$. For $j\in\{1,\ldots,q\}$, we define $\mu_j$ as the smallest index for which $\tau(e_{\mu_j})\geqslant 2\pi j/q$; we also set $\mu_0=1$. The desired partition consists of the $q$ sets some of which may be empty, and the $j$th of these sets consists of the vertices that lie between $e_{\mu_{j-1}}$ and $e_{\mu_j}$ again with respect to our enumeration.
\end{proof}

Observation~\ref{lemextrthin} shows that, for constant $\alpha$, a worst-case upper bound on the extension complexity of an $\alpha$-thin sequence gives an upper bound on $\wcc(n)$ which is at most a constant multiple worse than the bound in the $\alpha$-thin case. So we switch our attention to thin sequences to return to general polygons in Section~\ref{seccompl}.

\begin{defn}\label{defmanythin}
Let $v=(v_1,\ldots,v_n)$ be a thin cw-sequence. In what follows, the straight lines collinear to $v_1\wedge v_n$ are called \textit{horizontal}, and those orthogonal to $v_1\wedge v_n$ are \textit{vertical}. A non-vertical vector $\overrightarrow{a}$ is called \textit{left} if its projection onto $v_1\wedge v_n$ has the same direction as $\overrightarrow{v_1 v_n}$, and otherwise such a vector $\overrightarrow{a}$ is called \textit{right}. Also, we consider the line $v_1\wedge v_n$ which separates $\R^2$ into two half-planes, and we call the part containing $v$ the \textit{lower} half-plane, and the opposite one is called the \textit{upper} half-plane. Similarly to non-vertical vectors, we split the non-horizontal vectors into two classes, which are to be called the \textit{upgoing} and \textit{downgoing vectors}. A point $v$ is said to lie \textit{higher} than $u$ if the vector $\overrightarrow{uv}$ goes up. The quantity $\xc(\operatorname{conv} v)$ is also simply called the \textit{extension complexity} of $v$.
\end{defn}
 
\begin{defn}\label{defmanythin2}
Let $v=(v_1,\ldots,v_n)$ be a thin cw-sequence. For indexes $i,\hat{\imath},k, \hat{\jmath}, j$ satisfying $1\leqslant i<\hat{\imath}<k<\hat{\jmath}<j\leqslant n$, we define $\rho(v,i,\hat{\imath},k, \hat{\jmath}, j)$ as the ray passing from the point $(v_i\wedge v_{\hat{\imath}})\wedge(v_{\hat{\jmath}}\wedge v_j)$ towards $v_k$. 
 \end{defn}
 
\begin{figure}[h]
\includegraphics[width=12cm]{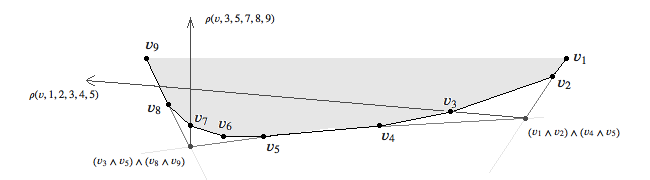}
\caption{
In this example, the ray $\rho(v,1,2,3,4,5)$ is left upgoing, and $\rho(v,3,5,7,8,9)$ is vertical upgoing.}
\label{figthin}
\end{figure}

\begin{lem}\label{lemleave}
Assume $v=(v_1,\ldots,v_n)$ is a thin sequence, and indexes $i,\hat{\imath},k, \hat{\jmath}, j$ are as in Definition~\ref{defmanythin2}. Then $\rho(v,i,\hat{\imath},k, \hat{\jmath}, j)\cap \operatorname{conv} v$ is a line segment $s$ in which $v_k$ is one of the endpoints, and $v_k$ lies between the apex of $\rho$ and the other endpoint of $s$.
\end{lem}

\begin{proof}
Let $\Phi$ be the set of all edge-defining inequalities of $\operatorname{conv} v$, and let $\rho'$ be the part of the ray $\rho$ obtained by removing all the points between $v_k$ and its apex $o=(v_i\wedge v_{\hat{\imath}})\wedge(v_{\hat{\jmath}}\wedge v_j)$. We write $\varphi_{1}\geqslant 0$, $\varphi_{2}\geqslant 0$ for the inequalities corresponding to the edges passing through $v_k$, and, since $v$ is thin, these inequalities should fail at $o$. So the quantities $\varphi_1$, $\varphi_2$ increase as we move from $o$ towards $v_k$, they vanish at $v_k$, and hence they should be positive on $\rho'$. All the other inequalities in $\Phi$ are strict at $v_k$, so a small part of $\rho'$ near its beginning lies in the interior of $\operatorname{conv} v$.
\end{proof}

\begin{defn}
We say that the ray $\rho$ as in Lemma~\ref{lemleave} \textit{enters} $\conv v$ at $v_k$ and \textit{leaves} $\conv v$ \textit{through} the other endpoint of $s$.
\end{defn}

\begin{ex}
The ray $\rho(v,1,2,3,4,5)$ on Figure~\ref{figthin} enters $\conv v$ at $v_3$ and leaves through a point on the segment $\conv\{v_8, v_9\}$. The ray $\rho(v,3,5,7,8,9)$ enters at $v_7$ and leaves through a point on $\conv\{v_1, v_9\}$.
\end{ex}


We finalize the section with a general property of acute diagrams whose bases have thin vertex sequences. We say that the \textit{continuation} of an oriented segment $\overrightarrow{uv}$ is the ray with the apex at $u$ and direction towards $v$.

\begin{lem}\label{acuteconv}
Let $\Delta$ be an acute diagram whose base has the vertices that form a thin sequence $v=(v_1,\ldots, v_n)$. We consider the orientation of $\Delta$ defined as in Remark~\ref{remrorient} with $v_1$ taken as a sink node. Consider an arbitrary inner node $u$ of $\Delta$ reachable by the oriented paths from the set $\Lambda(u)$ of base nodes, and assume that $v_n$ is not in $\Lambda(u)$. Let $\ell$ be a horizontal line that lies in the upper half-plane with respect to $v$. Assume that the continuation of any main arc coming from a vertex in $\Lambda(u)$ intersects $\ell$ at some point. If the set of all such intersection points has convex hull $\sigma$, then the continuation of the arc coming from $u$ intersects $\ell$ at a point in $\sigma$.
\end{lem}

\begin{proof}
Using the argument as in Lemma~\ref{lemmainint0}, we can find a base arc $\conv\{v,w\}$ with $v, w$ in $\Lambda(u)$ such that the main arcs passing from $v, w$ have the common end. With an inductive application of the procedure as in Lemma~\ref{lemlift}, we can assume without loss of generality that only main arcs can have $u$ as an endpoint, and also we can assume that there are exactly two such arcs $e_1, e_2$. In this case, it is straightforward to see that the arc coming from $u$ lies in the cone with the apex $u$ and extreme rays pointing towards the directions of $e_1, e_2$.
\end{proof}

\section{Good sequences}\label{secgoodseq}

We give several further auxiliary definitions to describe a situation that we call \textit{good} because the corresponding sequences allow small extended formulations.

\begin{defn}\label{defscat}
If $t, n$ are positive integers and $G$ is a subset of $\{1,\ldots,n\}$, then $G$ is called $t$\textit{-scattered} if for $g,\,g'\in G$ we have either $g=g'$ or $|g-g'|\geqslant t$.
\end{defn}

\begin{defn}\label{defgood0}
Assume $n\geqslant 5$ is an integer and $G\subset\{3,4,\ldots,n-2\}$ is a $3$-scattered subset. If $v=(v_1,\ldots,v_n)$ is a thin sequence, then the $G$\textit{-envelope} of $v$ is the sequence $v_G$ obtained from $v$ by replacing the points $v_{g-1}, v_g, v_{g+1}$ with $$o_g:=(v_{g-2}\wedge v_{g-1})\wedge (v_{g+1}\wedge v_{g+2}),$$ for all $g\in G$.
\end{defn}

As we see, the sequence $v_G$ is obtained from $v$ by a removal of three vertices and an addition of one new vertex, for all $g\in G$, which means that $v_G$ contains $n-2|G|$ vertices. More precisely, the polygon $\conv v_G$ is cut by all the edge-defining inequalities of $\operatorname{conv} v$ except $2|G|$ of them, namely, those corresponding to the edges $\operatorname{conv}\{v_{g-1},v_g\}$ and $\operatorname{conv}\{v_{g},v_{g+1}\}$ for $g\in G$. Further, we can see that $v_G$ is correct.

\begin{figure}[h]
\includegraphics[width=12cm]{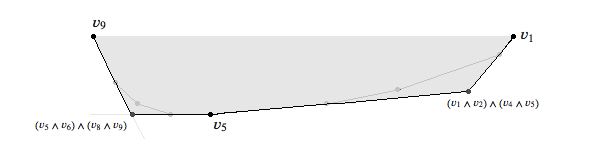}
\caption{The $G$-envelope of the sequence on Figure~\ref{figthin} with  $G=\{3,7\}$.}
\label{fig-g-good}
\end{figure}

\begin{obs}
If the sequence $v$ as in Definition~\ref{defgood0} is cw-ordered, then $v_G$ is a cw-ordered sequence with $n-2|G|$ vertices.
\end{obs}

\begin{proof}
Since $v$ is thin, the position of $o_g$ is such that $(v_{g-2}, o_g, v_{g+2})$ is a cw-ordered sequence, and the general statement follows by the induction on $|G|$.
\end{proof}


The motivation of the following definition lies in Lemma~\ref{lemgoodgood} below, which is a reformulation of Theorem~\ref{lemacutemult} in terms that are more comfortable to work with later.

\begin{defn}\label{defgood}
Assume $n\geqslant 5$ is an integer and $G\subset\{3,4,\ldots,n-2\}$ is a $3$-scattered subset. A correct sequence $v=(v_1,\ldots,v_n)$ is called $G$\textit{-good} if $v$ is thin and, additionally, there is an acute diagram with the base $\operatorname{conv} v_G$ such that, for any $g\in G$, the main edge passing from $o_g$ contains $v_g$.
\end{defn}

\begin{figure}[h]
\includegraphics[width=12cm]{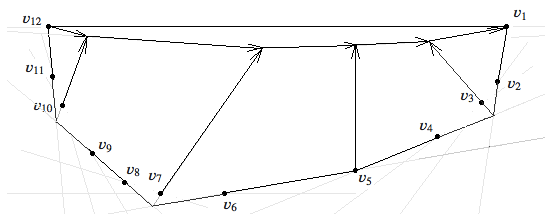}
\caption{An instance of Lemma~\ref{lemsameedge} with $n=12$ and $G=\{3,7,10\}$.}
\label{figout}
\end{figure}

A particular instance of Definition~\ref{defgood} is Figure~\ref{figout}, which corresponds to $n=12$ and shows a G-good sequence $(v_1,\ldots,v_{12})$ with $G=\{3,7,10\}$.

\begin{lem}\label{lemgoodgood}
Let $n\geqslant 5$ and $\delta\geqslant 1$ be integers, and let $v=(v_{11}, v_{21},\ldots,v_{n1})$ be a thin sequence. Let $G\subset\{3,4,\ldots,n-2\}$ be a $3$-scattered subset. Let $\mathcal{V}=(v_{gj})$ be an array of points in which the indexes $(g,j)$ run over $G\times\{1,\ldots,\delta\}$. We assume that, for any $j\in\{1,\ldots,\delta\}$, we obtain a good sequence if we replace $v_{g1}$ with $v_{gj}$ whenever $g\in G$. Then $\xc(\operatorname{conv} v\cup \mathcal{V})\leqslant n-|G|+\delta$.
\end{lem}

\begin{proof}
The idea of this proof is to apply Theorem~\ref{lemacutemult}; we define the polygon $P$ as in that theorem to be the convex hull of the $G$-envelope of $v$. The set $S$ as in Theorem~\ref{lemacutemult} is chosen as $\{o_g: g\in G\}$, where $$o_g:=(v_{g-2,1}\wedge v_{g-1,1})\wedge (v_{g+1,1}\wedge v_{g+2,1})$$
are the new vertices in the $G$-envelope of $v$. We pick $v_{g-1,1}$ and $v_{g+1,1}$ as the choice of $s', s''$ corresponding to the vertex $s=o_g$ as in Theorem~\ref{lemacutemult}. Since by the formulation of the lemma, a replacement of $v_{g1}$ with $v_{gj}$ returns a good sequence, and hence, in particular, a correct sequence, we see that the vertices $v_{g1},\ldots,v_{g\delta}$ lie in the interior of the triangle $\operatorname{conv} \{o_g, v_{g-1,1}, v_{g+1,1}\}$; these triangles are disjoint for different $g\in G$ again because the sequence $v$ is correct. The penultimate sentence of Theorem~\ref{lemacutemult} means in our setting that the sequence obtained by replacing $v_{g1}$ with $v_{gj}$ is $G$-good, and hence it is also valid by the assumptions of the current lemma.  So we have checked the conditions of Theorem~\ref{lemacutemult}, and we apply it to conclude that $\xc(\operatorname{conv} v\cup \mathcal{V})\leqslant (n-2|G|)+|G|+\delta=n-|G|+\delta$.
\end{proof}

We finalize this section with an important example of a good sequence.

\begin{lem}\label{lemsameedge}
Assume $n\geqslant 5$ is an integer and $G\subset\{3,4,\ldots,n-2\}$ is a $3$-scattered subset.  Assume $v=(v_1,\ldots,v_n)$ is a thin sequence such that, for all $g\in G$, the ray $\rho(v, {g-2}, g-1,g,g+1,g+2)$ as in Definition~\ref{defmanythin2} leaves $\operatorname{conv} v$ through the relative interior of the edge $\operatorname{conv}\{v_1, v_n\}$. Then $v$ is $G$-good.
\end{lem}

\begin{proof}
We recall that the $G$-envelope of $v$ is obtained by adding the point $$o_g:=(v_{g-2}\wedge v_{g-1})\wedge (v_{g+1}\wedge v_{g+2})$$ and removing $v_{g-1}, v_g, v_{g+1}$ for every $g\in G$. According to Definition~\ref{defgood}, we need to construct an acute diagram $\Delta$ in which the base is the $G$-envelope of $u$, and the main edge of $\Delta$ passing from $o_g$ contains $v_g$ for all $g\in G$.  We use Corollary~\ref{coracute} to build an acute diagram $\Delta$ with such a base in which

(1) for any $g\in G$, the main edge from any $o_g$ goes towards $v_g$,

(2) for any $k\notin\{1,n,g-1,g,g+1\}$, the main edge from $v_k$ passes in the vertical direction,

(3) the main edge from $v_n$ passes towards a point $u_0$ that lies sufficiently close to the middle of the edge $\operatorname{conv}\{v_1, v_n\}$ in the interior of $\operatorname{conv} v$.

We need to show that the main edge going from $o_g$ does actually contain $v_g$. We consider the orientation of $\Delta$ as in Remark~\ref{remrorient} with $v_1$ as a sink vertex. We proceed with a proof by contradiction; we assume that the desired statement is false, that is, there is an index $j\in G$ for which the main edge passing from $o_j$ does not contain $v_j$. This means that some edge $e$ of $\Delta$ has the ending point in the interior of $\operatorname{conv}\{o_g, v_j\}$, and hence the ending point of $e$ lies outside $\conv v$.

\textit{Case 1.} There exists an oriented path from $v_n$ to $e$. This is impossible because the oriented path starting at $v_n$ reaches $v_1$ without leaving $\conv v$ (in fact, it does not leave a sufficiently small neighborhood of $\conv\{v_1, v_n\}$ by the choice of $u_0$).

\textit{Case 2.} There is no oriented path from $v_n$ to $e$. According to Lemma~\ref{acuteconv}, the edges not having $v_n$ as a predecessor can exit $\conv v$ through the interior of $\conv\{v_1, v_n\}$ only, so any oriented path connecting a base vertex of $\Delta$ to the endpoint of $e$ should cross the interior of $\conv\{v_1, v_n\}$, and hence it should cross the oriented path between $v_n$ and $v_1$. Therefore, the situation falls into the already refuted Case~1; since these cases cover all possibilities, the proof is complete.
\end{proof}



\section{Extracting a slanted subsequence}\label{secextr}

We are ready to introduce the class of so-called \textit{slanted} sequences, which will be shown in Section~\ref{secslant} to contain sufficiently large subsequences of small extension complexity. In this section, we are going to explain how to extract a large slanted subsequence from a given sequence of large extension complexity. These two results will allow us to conclude the proof of Theorem~\ref{thrmain} in Section~\ref{seccompl}. 


\begin{defn}\label{defslant}
Let $v=(v_1,\ldots,v_t)$ be a thin cw-sequence, $\beta\in(0,\pi/2)$, $\delta\geqslant 0$. We say that $v$ is \textit{cw-slanted to} an angle $\beta$ with \textit{tolerance} $\delta$ if, for all $i,\hat{\imath},\hat{\jmath}, j$ satisfying $1\leqslant i<\hat{\imath}<\hat{\jmath}<j\leqslant t$, the ray ${\rho}(v,i,\hat{\imath},k,\hat{\jmath}, j)$ as in Definition~\ref{defmanythin2} satisfies
\begin{equation}\label{eqslant11}
\angle\left(\overrightarrow{v_{\hat{\jmath}}\wedge v_j},\,\,\rho(v,i,\hat{\imath},k,\hat{\jmath}, j)\right)<\beta
\end{equation}
for all $k$ satisfying $\hat{\imath}<k<\hat{\jmath}$ except for at most $\delta$ such values of $k$. 
\end{defn}

\begin{obs}\label{obsslant1}
Let $k'$, $k''$ be two indexes satisfying $\hat{\imath}<k'<k''<\hat{\jmath}$ in the notation of Definition~\ref{defslant}. If~\eqref{eqslant11} fails for $k=k''$, then it fails for $k=k'$ as well.
\end{obs}

\begin{proof}
Follows because $v$ is correct.
\end{proof}

\begin{obs}\label{obsslant2}
Let $\beta$, $\delta$, $v$ be as in Definition~\ref{defslant}. Then any ($\delta$+1)-scattered subsequence of $v$ is cw-slanted to the angle $\beta$ with tolerance $0$.
\end{obs}

\begin{proof}
Follows from Observation~\ref{obsslant1}.
\end{proof}

\begin{remr}\label{defslant2}
If a correct sequence $v=(v_1,\ldots,v_t)$ is not cw-ordered, which means that it is ccw-ordered, then we say that $v$ is cw-slanted with parameters as in Definition~\ref{defslant} if the reversed sequence satisfies the assumptions of this definition.
\end{remr}

\begin{remr}
We define the dual concept of \textit{ccw-slantedness} by replacing every occurrence of `\textit{cw}' with `\textit{ccw}' and vice versa in Definition~\ref{defslant} and Remark~\ref{defslant2}.
\end{remr}


The goal of the remaining part of this section is to prove Lemma~\ref{lemextrrob} below; we need to begin with an important special case of it.

\begin{lem}\label{lemauxextr}
Let $G\subset\{3,4,\ldots,n-3\}$ be a $4$-scattered subset with $n\geqslant 5$, and let $v=(v_1,\ldots,v_n)$ be a thin cw-sequence. Assume that, for any $g\in G$, the ray $\rho(v,g-2,g-1,g,g+1,g+2)$ is left, and $v_g$ is located to the right of the point 
$$\omega_g:=(v_{g-2}\wedge v_{g-1})\wedge (v_{g+2}\wedge v_{g+3}).$$ 
Then there exist arbitrarily high points $u_1$ and $u_n$ such that $u_1$ lies on the ray passing from $v_2$ towards $v_1$, and $u_n$ lies on the ray passing from $v_{n-1}$ towards $v_n$, for which the sequence $(u_1, v_2,\ldots,v_{n-1}, u_n)$ is $G$-good.
\end{lem}

\begin{figure}[h]
\includegraphics[width=12cm]{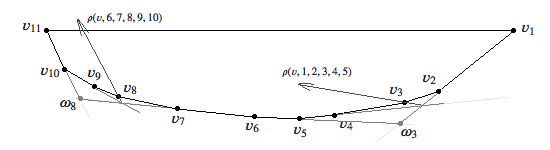}
\caption{An instance of Lemma~\ref{lemauxextr} with $n=11$ and $G=\{3,8\}$.}
\label{figextr11}
\end{figure}

\begin{proof}
For any $g\in G$, we set 
$$o_g:=(v_{g-2}\wedge v_{g-1})\wedge (v_{g+1}\wedge v_{g+2})$$
and denote by $\ell_g$ the vertical line passing through $\omega_g$. Also, we define
$$\alpha_g:=\ell_g\wedge (o_g\wedge v_g),$$
and we set the points $u_1, u_n$ as in the formulation of the lemma so that the segment connecting them is horizontal and lies higher than the point $\alpha_g$ for any $g\in G$. Proving that $(u_1, v_2,\ldots,v_{n-1}, u_n)$ is a $G$-good sequence requires an acute diagram whose base is the convex hull of the $G$-envelope of $(u_1, v_2,\ldots, v_{n-1}, u_n)$, so the base vertices of such a diagram can be separated into the following three types:

(A) $u_1$, $u_n$,

(B) $o_g$ and $v_{g+2}$ with $g\in G$,

(C) all other vertices.

\noindent By the definition of the G-envelope, we note that the type C vertices should include those $v_j$ for which $j\notin\{1,n\}$ and also $j\notin\{g-1,g,g+1,g+2\}$ for any $g\in G$.

Now we are ready to construct the diagram $\Delta$ with the use of Lemma~\ref{coracute}. The main edge passing from $u_1$ gets assigned the direction towards an inner point close to the middle of the segment $\conv\{u_1, u_n\}$; the vertices of type C are assigned the vertical direction of the corresponding main edges. Finally, the vertices of type B have main edges passing towards $\alpha_g$ with $g\in G$.

\begin{figure}[h]
\includegraphics[width=12cm]{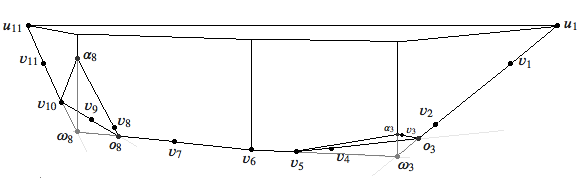}
\caption{The proof of Lemma~\ref{lemauxextr} with the polygon on Figure~\ref{figextr11}.}
\label{figextr12}
\end{figure}

In order to complete the proof, we need to check that the main edge passing from $o_g$ does indeed contain $v_g$ for any $g\in G$. Indeed, a stronger fact follows by the induction on the cardinality of $G$, namely, that the two main edges passing from $o_g$ and $v_{g+2}$ meet at $\alpha_g$, and the third arc adjacent to $\alpha_g$ goes vertically until it meets the path between $u_1$ and $u_n$ close to the segment $\conv\{u_1, u_n\}$, and, similarly, a main edge passing from a vertex of type C goes in the vertical direction until it meets the path between $u_1$ and $u_n$, as it is the case with $v_6$ on Figure~\ref{figextr12}.
\end{proof}

Now we are ready to proceed with the main result of this section.

\begin{lem}\label{lemextrrob}
Let $\delta, \tau, m$ be positive integers and $n=8\tau m$. Let $\alpha,\beta$ be positive reals with $\pi/2>\beta\geqslant 2\alpha$. Assume $v=(v_1,\ldots,v_n)$ is an $\alpha$-thin sequence such that any its subsequence with $m$ points is neither cw-slanted nor ccw-slanted relative to the angle $\beta$ and tolerance $2\delta$. Then $v$ has a subsequence with $(6+\delta)\tau$ points with extension complexity not exceeding $6\tau+\delta+1$.
\end{lem}

\begin{proof}
We assume without loss of generality that the vertices of $v$ are enumerated in the cw-order. The assumption on the lack of slanted subsequences implies that
$$(v_1,\ldots, v_m),\,\,(v_{2m+1},\ldots,v_{3m}),\,\,\ldots,\,\,(v_{(8\tau-2)m+1},\ldots,v_{(8\tau-1)m})$$
are not ccw-slanted with respect to the angle $\beta$ and tolerance $2\delta$, and also
$$(v_{m+1},\ldots, v_{2m}),\,\,(v_{3m+1},\ldots,v_{4m}),\,\,\ldots,\,\,(v_{(8\tau-1)m+1},\ldots,v_{8\tau m})$$
are not cw-slanted relative to the same parameters. This means that, for any $q\in\{0,\ldots,4\tau-1\}$, we can find $i_q$, $\hat{\imath}_q$, $\hat{\jmath}_q$, $j_q$ and $i'_q$, $\hat{\imath}'_q$, $\hat{\jmath}'_q$, $j_q$ satisfying
$$2qm+1\leqslant i_q<\hat{\imath}_q<\hat{\jmath}_q<j_q\leqslant (2q+1)m,$$
$$(2q+1)m+1\leqslant i'_q<\hat{\imath}'_q<\hat{\jmath}'_q<j'_q\leqslant (2q+2)m,$$
for which the ray $\rho_{k_q}:=\rho(i_q,\hat{\imath}_q,k_q,\hat{\jmath}_q,j_q)$ satisfies
\begin{equation}\label{eqray1}
\angle\left(\overrightarrow{v_{\hat{\imath}_q} v_{i_q}},\rho_{k_q}\right)\geqslant\beta
\end{equation}
for all $k_q$ in a set $\Delta_q\subseteq\{\hat{\imath}_q+1,\ldots,\hat{\jmath}_q-1\}$ of cardinality at least $2\delta$, and the ray $\rho_{k'_q}:=\rho(i'_q,\hat{\imath}'_q,k'_q,\hat{\jmath}'_q,j'_q)$ satisfies
$$\angle\left(\overrightarrow{v_{\hat{\jmath}'_q} v_{j'_q}},\rho_{k'_q}\right)\geqslant\beta$$
for all $k'_q$ in a set $\Delta'_q\subseteq\{\hat{\imath}'_q+1,\ldots,\hat{\jmath}'_q-1\}$ of cardinality at least $2\delta$.
Since $v$ is cw-ordered, every point in $\Delta_q$ lies to the right of any point of $\Delta'_q$, so for the point
$$\omega_q:=(v_{i_q}\wedge v_{\hat{\imath}_q})\wedge(v_{\hat{\jmath}'_q}\wedge v_{j'_q})$$
we have at least one of the two options:

(a) every point of $\Delta_q$ lies to the right of $\omega_q$,

(b) every point of $\Delta'_q$ lies to the left of $\omega_q$.

We can assume without loss of generality that the quantity of those $q$ for which (a) holds is not less than the corresponding quantity for (b), because otherwise this property can be satisfied on the sequence obtained from $v$ by the axial symmetry. So we can find a set $Q\subset\{0,\ldots,4\tau-1\}$ of cardinality $2\tau$ such that any $q\in Q$ has the property (a). Further, an index $q\in Q$ is called \textit{left-admissible} if there is a set $\mathcal{L}_q\subset\Delta_q$ of cardinality $\delta$ consisting of those $k_q\in\Delta_q$ for which the ray $\rho_{k_q}$ is left. Otherwise, an index $q\in Q$ is called \textit{right-admissible}. The situation splits into two possible cases, which we need to consider separately.

\textit{Case 1.}  Assume that the quantity of the left-admissible indexes in $Q$ is not less than that of the right-admissible indexes. This means that we can pick a set $Q_l$ of cardinality $\tau$ such that, for any $q$ in $Q_l$, there exists a set $\mathcal{L}_q\subset\Delta_q$ of $\delta$ values of $k_q$ for which the ray $\rho_{k_q}$ is left. Then we define the subsequence $u$ of $v$ containing the points
$$
\bigcup\limits_{q\in Q_l}\left(\L_q\cup\left\{v_{i_q}, v_{\hat{\imath}_q}, v_{\hat{\jmath}_q}, \alpha_q, v_{j'_q}\right\}\right),
$$
where $$\alpha_q:=(v_{\hat{\jmath}_q}\wedge v_{j_q})\wedge(v_{\hat{\jmath}'_q}\wedge v_{j'_q}).$$
For any choice of $k_q\in\mathcal{L}_q$, the subsequence $u_q$ of $u$ consisting of
$$\bigcup\limits_{q\in Q_l}\left\{v_{i_q}, v_{\hat{\imath}_q}, v_{k_q}, v_{\hat{\jmath}_q}, \alpha_q, v_{j'_q}\right\}$$
satisfies the assumptions of Lemma~\ref{lemauxextr} with $G$ being the set of the positions of the $v_{k_q}$'s. According to this lemma, we see that the sequence $u_q$ becomes $G$-good after \textit{stretching the first and last vertices}, by which we mean that one replaces these vertices by two higher points, one of which is taken on the ray passing from the second point towards the first point, and the other one is chosen on the ray passing from the penultimate point towards the last point, respectively. Also, let $u'$ be a sequence obtained from $u$ by stretching the first and last vertices; this sequence has $(5+\delta)\tau$ vertices. According to Lemma~\ref{lemgoodgood}, we can get
$$\xc(\conv u')\leqslant 6\tau - \tau + \delta=5\tau+\delta.$$
Now we can cut $\conv u'$ by the line joining the first and last points of $u$, which acts like a takeback of the stretching operation. In other words, this leaves us with $\conv u$, which we further cut along $v_{j_q}\wedge v_{\hat{\jmath}_q'}$ for all $q\in Q_l$, thus adding the $v_{j_q}, v_{\hat{\jmath}_q'}$ to $u$ and getting rid of all the points $\alpha_q$ defined above. So we get from $u'$ a subsequence of $v$ of cardinality $(6+\delta)\tau$ for the cost of $\tau+1$ additional cuts, each of which can worsen our bound only by at most one, according to Observation~\ref{lemeasy1}, so this subsequence has extension complexity not exceeding $6\tau+\delta+1$.

\textit{Case 2.} The set $Q$ has at least as many right-admissible indexes as left-admissible indexes. Similarly to the consideration in the first paragraph of Case~1, we find a set $Q_r$ of cardinality $\tau$ 
such that, for any $q$ in $Q_r$, there exists a set $\mathcal{R}_q\subset\Delta_q$ of $\delta$ values of $k_q$ for which the ray $\rho_{k_q}$ is not left. We define the subsequence $w$ of $v$ containing the points
$$
\bigcup\limits_{q\in Q_r}\left(\Rr_q\cup\left\{v_{i_q}, v_{\hat{\imath}_q}, v_{\hat{\jmath}_q}, v_{j_q}\right\}\right),
$$
and we denote by $w_r, w_l$ the first and last points of $w$, respectively. According to Observation~\ref{obsaddang}, we have
$$\angle\left(\rho_{k_q},\overrightarrow{w_rw_l}\right)\geqslant\angle\left(\rho_{k_q},\overrightarrow{v_{\hat{\imath}_q} v_{i_q}}\right)-\angle\left(\overrightarrow{v_{\hat{\imath}_q} v_{i_q}},\overrightarrow{w_rw_l}\right),$$
and using the bound of~\eqref{eqray1} for the first summand and the $\alpha$-thinness for the second summand, we get
$$\angle\left(\rho_{k_q},\overrightarrow{w_rw_l}\right)>\beta-\alpha\geqslant\alpha,$$
so for $k_q\in\mathcal{R}_q$, the ray $\rho_{k_q}$ has to leave $\conv w$ through the interior of $\conv\{w_l, w_r\}$. According to Lemma~\ref{lemsameedge}, the sequence
$$
\bigcup\limits_{q\in Q_r}\left\{v_{i_q}, v_{\hat{\imath}_q}, v_{k_q}, v_{\hat{\jmath}_q}, v_{j_q}\right\}
$$
is $G$-good for any choice of $k_q\in\mathcal{R}_q$, where $G$ is the set of the positions of the $v_{k_q}$'s. Using Lemma~\ref{lemgoodgood}, we conclude that the extension complexity of $w$ is at most $5\tau-\tau+\delta=4\tau+\delta.$ Since the cardinality of $w$ is $(4+\delta)\tau$, we can add $2\tau$ arbitrary points of $v$ to $w$, and the resulting sequence will have extension complexity not exceeding $6\tau+\delta$ according to Lemma~\ref{lemeasy2}. 
\end{proof}

\section{Decreasing slanted sequences}\label{secdecr}

In the previous section, we introduced the notion of a slanted sequence. Roughly speaking, we showed that a polygon of large extension complexity should admit a sufficiently large slanted subsequence. This result allows us to concentrate on slanted sequences, and the following definition turns out to be useful. We recall that the notation $\dist(A, B)$ stands for the distance between two sets $A, B\subset\R^2$.

\begin{defn}\label{defdecr}
Let $u=(u_1, u_1', w_{12}, u_2, u_2', w_{23},\ldots, w_{s-1, s}, u_s, u_s')$ be a thin cw-sequence with $3s-1$ points. We say that $u$ is \textit{cw-decreasing} if the inequality
$$\dist(w_{q-1, q}, u_q\wedge u_q')>\dist(w_{q, q+1},u_q\wedge u_q')$$
holds for any $q\in\{2,3,\ldots,s\}$.
\end{defn}

\begin{figure}[h]
\includegraphics[width=8cm]{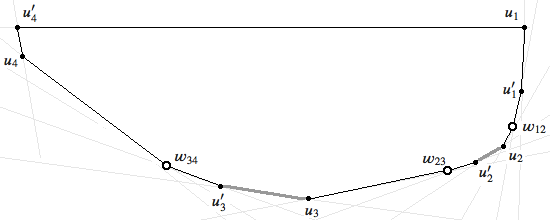}
\caption{An instance of Definition~\ref{lemsameedge} with $s=4$: Among two consecutive punctured points $w_{q-1\, q}$ and $w_{q\, q+1}$, the left one is closer to the straight line $u_q\wedge u_q'$ containing the thick edge.}
\label{figdecr}
\end{figure}

Now we are going to show that cw-decreasing sequences admit acute extensions under an additional mild assumption on the positions of the vertices. The proof of the following lemma uses an idea similar to that of Lemma~\ref{lemsameedge}. 

\begin{lem}\label{lemdecgood}
Let $u=(u_1, u_1', w_{12}, \ldots, w_{s-1, s}, u_s, u_s')$ be a thin cw-decreasing sequence written in the cw-order. Assume that the ray passing from
$$\omega_{q, q+1}:=(u_q\wedge u_q')\wedge(u_{q+1}\wedge u_{q+1}')$$
towards $w_{q, q+1}$ is left, for any $q\in\{1,\ldots,s-1\}$. Then $u$ is $G$-good, where $G=\{3,6,\ldots,3s-3\}$ is the set of the positions of the $w_{q, q+1}$'s.
\end{lem}

\begin{proof}
To prove the desired result, we need to construct an acute diagram $\Delta$ whose base is the $G$-envelope of $u$, and the main edge of $\Delta$ passing from every $\omega_{q, q+1}$ should contain $w_{q, q+1}$. According to Definition~\ref{defgood}, the $G$-envelope of $u$ is $(u_1,w_{12}, w_{23},\ldots,w_{s-1, s},u_s')$, and, as in the proof of Lemma~\ref{lemsameedge}, we use Corollary~\ref{coracute} to build an acute diagram $\Delta$ in which

(1) for any $q\in\{1,\ldots,s-1\}$, the main edge from any $\omega_{q,q+1}$ goes towards $w_{q,q+1}$,

(2) the main edge from $u_s'$ passes towards a point $u_0$ that lies sufficiently close to the middle of the edge $\operatorname{conv}\{u_1, u_s'\}$ in the interior of $\operatorname{conv} u$.

We need to show that the main edge going from $\omega_{q, q+1}$ does actually contain $w_{q, q+1}$. We consider the orientation of $\Delta$ as in Remark~\ref{remrorient} with $u_1$ as a sink vertex. We proceed with a proof by contradiction; we assume that the desired statement is false, that is, there is an index $j\in\{1,\ldots,s-1\}$ for which the main edge passing from $\omega_{j,j+1}$ does not contain $w_{j,j+1}$. This means that some edge $e$ of $\Delta$ has the ending point $\varepsilon$ in the interior of $\operatorname{conv}\{\omega_{j,j+1},w_{j,j+1}\}$, and hence this point $\varepsilon$ lies outside $\conv u$. Furthermore, we assume that $j$ is the minimal index that possesses this property.

\textit{Case 1.} There exists an oriented path from $u_s'$ to $e$. Similarly to the proof of Lemma~\ref{lemsameedge}, this is impossible because the oriented path from $u_s'$ to $u_1$ lies sufficiently close to $\conv\{u_1, u_s'\}$.

\textit{Case 2.} There exists an oriented path from some $\omega_{t\, t+1}$ to $e$ with $t>j$. This is again impossible because $w_{j\,j+1}$ is located to the right of $\omega_{t\,t+1}$, and every edge of $\Delta$ that is not reachable from $u_s'$ is left because of Lemma~\ref{acuteconv}.

Since Cases~1 and~2 are invalid, we can use Lemma~\ref{lemextra} and conclude that $\varepsilon$ is reachable by an oriented path from $\omega_{j-1, j}$. We denote this path by $\pi$, and we conclude, by the minimality of $j$, that the main edge of $\omega_{j-1\,j}$ passes through $w_{j-1\,j}$, which means that $w_{j-1\,j}$ lies on $\pi$. 
Using the thinness of $u$,
we conclude that the distance from a point $x\in\pi$ to $u_j\wedge u_j'$ grows as we move $x$ towards the direction of $\pi$. Since $u$ is cw-decreasing, we have
$$\dist(w_{j,j+1}, u_j\wedge u_j')<\dist(w_{j-1,j}, u_j\wedge u_j'),$$
so the point $w_{j,j+1}$ occurs on $\pi$ earlier than $w_{j-1,j}$, which is a contradiction.
\end{proof}

The main result of this section is a corollary of Lemma~\ref{lemdecgood}.

\begin{cor}\label{corcwdecr}
Let $\delta,\tau$ be positive integers; let $\alpha,\beta$ be positive reals with $\alpha+\beta\leqslant \pi/2$. Suppose $n=4\delta\tau+2$ and let $v=(v_1,\ldots,v_n)$ be an $\alpha$-thin cw-sequence. If $v$ is cw-slanted to the angle $\beta$ with tolerance $2\delta$, and if the subsequence
$$v'=(v_1, v_2, v_{4\delta}, v_{4\delta+1}, v_{4\delta+2}, v_{8\delta}, \ldots, v_{n-1}, v_n)$$
is cw-decreasing, then $v$ possesses a subsequence of $2\delta\tau+2$ points with extension complexity not exceeding $2\tau+2\delta$.
\end{cor}

\begin{proof}
According to Definition~\ref{defslant}, every choice of the index $q=\{0,\ldots,\tau-1\}$ allows a set
$\Delta_q\subset\{4q\delta+3,\ldots,4q\delta+4\delta\}$
of cardinality $2\delta-2$ such that the ray
$$\rho_{k_q}:=\rho(v,4q\delta+1,{4q\delta+2},{k_q},{4(q+1)\delta+1},{4(q+1)\delta+2})$$
as in Definition~\ref{defmanythin2} satisfies 
$$\angle(\rho_{k_q},\overrightarrow{v_{4(q+1)\delta+1}v_{4(q+1)\delta+2}})<\beta$$ for $k_q\in\Delta_q$. According to
Observation~\ref{obsaddang}, the angle $\angle(\rho_{k_q},\overrightarrow{v_1v_n})$ does not exceed
$$\angle(\rho_{k_q},\overrightarrow{v_{4(q+1)\delta+1}v_{4(q+1)\delta+2}})+\angle(\overrightarrow{v_1v_n},\overrightarrow{v_{4(q+1)\delta+1}v_{4(q+1)\delta+2}})<\beta+\alpha\leqslant\pi/2,$$
so the ray $\rho_{k_q}$ is left for any $k_q\in\Delta_q$. We also have
\begin{equation}\label{eqeqeq1}
\dist(v_{k_q},v_{4(q+1)\delta+1}\wedge v_{4(q+1)\delta+2})\geqslant \dist(v_{4(q+1)\delta},v_{4(q+1)\delta+1}\wedge v_{4(q+1)\delta+2}),
\end{equation}
\begin{equation}\label{eqeqeq2}
\dist(v_{4(q+2)\delta},v_{4(q+1)\delta+1}\wedge v_{4(q+1)\delta+2})\geqslant \dist(v_{k_{q+1}},v_{4(q+1)\delta+1}\wedge v_{4(q+1)\delta+2})
\end{equation}
for all $k_q\in\Delta_q$ and $k_{q+1}\in\Delta_{q+1}$ because $v$ is correct, and also
$$
\dist(v_{4(q+1)\delta},v_{4(q+1)\delta+1}\wedge v_{4(q+1)\delta+2})
>\dist(v_{4(q+2)\delta},v_{4(q+1)\delta+1}\wedge v_{4(q+1)\delta+2})
$$
because $v'$ is cw-decreasing. Putting the last inequality together with~\eqref{eqeqeq1}--\eqref{eqeqeq2}, we get
$$
\dist(v_{k_q},v_{4(q+1)\delta+1}\wedge v_{4(q+1)\delta+2})
>\dist(v_{k_{q+1}},v_{4(q+1)\delta+1}\wedge v_{4(q+1)\delta+2}),
$$
which means that the sequence
$(v_1,v_2,v_{k_0},v_{4\delta+1}, v_{4\delta+2}, v_{k_1},\ldots,v_{n-1},v_n)$
is cw-decreasing for any choice of $v_{k_q}\in\Delta_q$, and by Lemma~\ref{lemdecgood} this sequence is $G$-good, where $G$ is the set of the positions of the $v_{k_q}$'s. It remains to apply Lemma~\ref{lemgoodgood} to see that the subsequence formed by the set
$$\{v_1,v_2,v_{4\delta+1}, v_{4\delta+2}, \ldots,v_{n-1},v_n\}\cup\Delta_0\ldots\cup\Delta_{\tau-1}$$
has extension complexity not exceeding $(3\tau+2)-\tau+(2\delta-2)=2\tau+2\delta$.
\end{proof}

\section{Extensions for slanted sequences}\label{secslant}

In this section, we finalize our discussion of slanted sequences by proving a general result concerning their extension complexities. The following auxiliary lemma is quite straightforward.

\begin{lem}\label{lemrayadd}
Let $u=(u_2,u_3,u_4,u_5,u_6)$ be a thin cw-ordered sequence. Assume that $\rho$ is a ray with apex at $u_4$. If $\rho$ intersects $\conv\{u_5, u_6\}$, then $$\angle\left(\overrightarrow{u_2 u_3}, \rho\right)<\angle\left(\overrightarrow{u_2 u_3}, \overrightarrow{u_5 u_6}\right).$$
\end{lem}

\begin{proof}
Let $w_6$ be the intersection point of $\conv\{u_5, w_6\}$ and $\rho$. We get
$$\angle\left(\overrightarrow{u_2 u_3}, \rho\right)=\angle u_2 u_3 u_4+\angle u_3u_4w_6,$$
$$\angle\left(\overrightarrow{u_2 u_3}, \overrightarrow{u_5 w_6}\right)=\angle u_2u_3u_4+\angle u_3 u_4 u_5+\angle u_4u_5w_6,$$
and we are done because $\angle u_3 u_4 u_5+\angle u_4u_5w_6>\angle u_3u_4w_6$.
\end{proof}

The following two lemmas are based on direct computations. In Lemma~\ref{lemrelation} below, it is useful for us to develop a representation of a thin sequence in Cartesian coordinates different from that of Definition~\ref{defmanythin}. Namely, for a given thin sequence $v=(v_2,\ldots,v_t)$ with $t\geqslant 4$ that is written in the cw-order, we can construct a congruence transformation of $\R^2$ that sends 

(i) $v_2$ to a point $(a,0)$ with $a>0$, $v_3$ to the point $(0,0)$,

(ii) $v_4,\ldots, v_t$ to points with positive second coordinates.

In fact, the images of two distinct points as in (i) determine the congruence transformation up to a symmetry, and since our sequence corresponds to a convex polygon, the remaining points should have second coordinates of the same sign.

\begin{figure}[h]
\includegraphics[width=8cm]{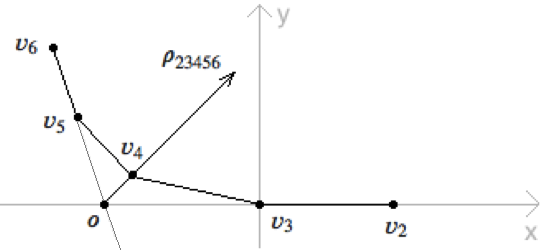}
\caption{The coordinate system as in the proof of Lemma~\ref{lemrelation}.}
\label{figcoord}
\end{figure}

\begin{lem}\label{lemrelation}
Let $v=(v_1,\ldots,v_t)$ be an $\alpha$-thin sequence cw-slanted to an angle $\beta$ with tolerance $0$, where $\alpha+\beta\leqslant\pi/2$. Then, for any tuple of indexes $i,\hat{\imath},k,\hat{\jmath}, j$ satisfying $1\leqslant i<\hat{\imath}<k<\hat{\jmath}<j\leqslant t$, we have
$$\dist(v_{\hat{\imath}},v_k)\sin(\angle v_i v_{\hat{\imath}} v_k)>\dist(v_k,v_{\hat{\jmath}})\sin(\angle v_k v_{\hat{\jmath}} v_j).$$
\end{lem}

\begin{proof}
We introduce the notation $L_{pq}=\dist(v_p,v_q)$, 
$a_{pqr}=\pi-\angle v_p v_q v_r$ and $a_s=a_{s-1,s,s+1}$, and also we write $\rho_{i\hat{\imath}k\hat{\jmath} j}=\rho(v,i,\hat{\imath},k,\hat{\jmath}, j)$. We assume without loss of generality that $v$ is written in the cw-order, and also that $i=2$, $\hat{\imath}=3$, $k=4$, $\hat{\jmath}=5$, $j=6$. The congruence transformation described above allows us to take
$$v_2=(L_{23},0),\,\,\,\,v_{3}=(0,0),\,\,\,\,v_4=(-L_{34}\cos a_{3}, L_{34}\sin a_{3}),$$
$$v_{5}=v_4+(-L_{45}\cos (a_{3}+a_{4}), L_{45}\sin (a_{3}+a_{4})),$$
$$v_{6}=v_{5}+(-L_{56}\cos (a_{3}+a_{4}+a_{5}), L_{56}\sin (a_{3}+a_{4}+a_{5})).$$
Now we want to compute the angle between the ray $\rho_{23456}$ and the vector $\overrightarrow{v_2v_3}$, which goes towards the negative direction of the x-axis. This ray has the apex at $o:=(v_2 \wedge v_{3})\wedge(v_{5}\wedge v_6)$ and its direction is towards $v_4$; a direct computation gives
$$o=\left(\frac{-L_{45}\sin a_{5} - L_{34}\sin(a_{4} + a_{5})}{\sin(a_{3} + a_{4} + a_{5})},0\right).$$
The $\tan$ of the angle between $\rho_{23456}$ and $\overrightarrow{v_2v_3}$ equals
\begin{equation}\label{oldereq}
\frac{L_{34} \sin a_{3} \sin(a_{3} + a_{4} + a_{5})}{L_{34} \sin a_{3} \cos(a_{3} + a_{4} + a_{5}) - L_{45}\sin a_{5}},
\end{equation}
and we are done because the denominator of~\eqref{oldereq} should be positive (if this was not the case, we would have $\angle(\rho_{23456},\overrightarrow{v_2v_3})\geqslant\pi/2$ and hence $\angle(\rho_{23456},\overrightarrow{v_5v_6})>\pi/2-\alpha\geqslant\beta$, which is a contradiction to the cw-slantedness). 
\end{proof}

We proceed with a deeper result based on computations similar to Lemma~\ref{lemrelation}.

\begin{lem}\label{lemout}
Let $v=(v_1,\ldots,v_{10})$ be a $\pi/6$-thin sequence cw-slanted relative to the angle $\pi/3$ with tolerance $0$. We assume that $\gamma_3\geqslant1$, $\gamma_4\geqslant1$, $\gamma_5\geqslant1$, $\gamma_6\geqslant1$, where 
\begin{equation}\label{eqgamma}
\gamma_\tau=\frac{\dist(v_{\tau+3}, v_{\tau+1}\wedge v_{\tau+2})}{\dist(v_{\tau}, v_{\tau+1}\wedge v_{\tau+2})}.
\end{equation}
Then there are indexes $i,\hat{\imath},k,\hat{\jmath}, j$ satisfying $2\leqslant i<\hat{\imath}<k<\hat{\jmath}<j\leqslant 9$ such that $\rho(v,i,\hat{\imath},k,\hat{\jmath}, j)$ leaves $\operatorname{conv} v$ through the relative interior of $\operatorname{conv}\{v_1, v_{10}\}$.
\end{lem}

\begin{proof}
We continue to use the notation of the proof of Lemma~\ref{lemrelation}. As we saw from the cw-slantedness, the ray $\rho_{23456}$ should be left, so the only possibility for it to fail the desired statement, that is, not to leave $\operatorname{conv} v$ through the relative interior of $\operatorname{conv}\{v_1, v_{10}\}$, is that $\rho_{23456}$ intersects $\operatorname{conv}\{v_6, v_{10}\}$.
Using Lemma~\ref{lemrayadd} with $v_2, v_3, v_4,v_6,v_{10}$ taking the roles of $u_2, u_3, u_4,  u_5, u_6$, respectively, and with the part of $\rho_{23456}$ beginning at $v_4$ in the role of $\rho$, we get
%
%
\begin{equation}\label{eqaux1122}
\angle(\rho_{23456},\overrightarrow{v_2v_3})<\angle(\overrightarrow{v_6v_{10}},\overrightarrow{v_2v_3}).
\end{equation}
We note that the right-hand side of~\eqref{eqaux1122} is  $a_3+a_4+a_5+a_{5\,6\,10}$, and the left-hand side of~\eqref{eqaux1122} was computed in the proof of Lemma~\ref{lemrelation} and equals the $\arctan$ of the expression~\eqref{oldereq}. We get
\begin{equation}\label{eqaux11}
\arctan\left(\frac{L_{34} \sin a_{3} \sin(a_{3} + a_{4} + a_{5})}{L_{34} \sin a_{3} \cos(a_{3} + a_{4} + a_{5}) - L_{45}\sin a_{5}}\right)-a_3-a_4-a_5<a_{5\,6\,10}.
\end{equation}
We have $L_{56}\sin a_{5\,6\,10}<L_{45}\sin a_4$ by the result of of Lemma~\ref{lemrelation}, so the right-hand side of~\eqref{eqaux11} can be replaced by $\arcsin(L_{45}\sin a_4/L_{56})$. We do this and further take the $\tan$ function of both sides of~\eqref{eqaux11}; we recall that $\tan(x+y)=(\tan x+\tan y)/(1-\tan x\tan y)$ and $\tan(\arcsin x)=x/\sqrt{1-x^2}$. We get 
$$\frac{L_{45}\sin a_5 \sin(a_3+a_4+a_5)}{L_{34}\sin a_3-L_{45}\sin a_5\cos (a_3+a_4+a_5)}<\frac{L_{45}\sin a_4}{\sqrt{L_{56}^2-L_{45}^2\sin^2 a_4}},$$
and now we are able to remove the second listed summand of the denominator of the left-hand side. Further elementary transformations yield
$$\sin a_5 \sin(a_3+a_4+a_5)\sqrt{L_{56}^2-L_{45}^2\sin^2 a_4}<L_{34}\sin a_3\sin a_4.$$
We note that $L_{45}\sin a_5<L_{34}\sin a_3$ by the result of Lemma~\ref{lemrelation}, which gives
\begin{equation}\label{eqnewnew11}
\sin a_5 \sin(a_3+a_4+a_5)\sqrt{L_{56}^2-\frac{L_{34}^2\sin^2 a_4\sin^2 a_3}{\sin^2 a_5}}<L_{34}\sin a_3 \sin a_4,
\end{equation}
and now we use the equation~\eqref{eqgamma} with $\tau=3$, which means that
$$L_{34}\sin a_4=\frac{L_{56}\sin a_5}{\gamma_3},$$
and we make the replacement of $L_{34}\sin a_4$ in~\eqref{eqnewnew11} to get
$$\sin a_5 \sin(a_3+a_4+a_5)\sqrt{L_{56}^2-\frac{L_{56}^2 \sin^2 a_3}{\gamma_3^2}}<\frac{L_{56}\sin a_3\sin a_5}{\gamma_3},$$
and further elementary cancellations give
\begin{equation}\label{qwertyu1}
\sin(a_3+a_4+a_5)\sqrt{\gamma_3^2-\sin^2 a_3}<\sin a_3.
\end{equation}
This inequality is obviously false if the square root is greater than or equal to one, so we have $\gamma_3<\sqrt{1+\sin^2 a_3}$, which implies 
\begin{equation}\label{eqeas1}
\gamma_3<1.5.
\end{equation}
Taking into account the inequality $\gamma_3\geqslant1$ from the assumptions of the lemma and the inequality  $a_3<\pi/6$ that comes from the $\pi/6$-thinness of $v$, we transform~\eqref{qwertyu1} into
$$\sin(a_3+a_4+a_5)\sqrt{1-(\sin \pi/6)^2}<\sin a_3,$$ from which it can be deduced that 
\begin{equation}\label{eqeas2}
a_4<a_3
\end{equation}
by standard tools of calculus.

Now we apply similar considerations with the rays $\rho_{56789}$ and $\rho_{23458}$ in the role of the ray $\rho_{23456}$ above. In the above consideration, the only case of the inequality $\gamma_i\geqslant 3$ that we used was $i=3$, so the ability to make the same computations for $\rho_{56789}$ is justified by the analogous inequality $\gamma_6\geqslant 1$ in the formulation. Similarly, the case of $\rho_{23458}$ is guaranteed by
$$\frac{\dist(v_{8}, v_{4}\wedge v_{5})}{\dist(v_{3}, v_{4}\wedge v_{5})}>\frac{\dist(v_{6}, v_{4}\wedge v_{5})}{\dist(v_{3}, v_{4}\wedge v_{5})}=\gamma_3\geqslant1.$$

In particular, if the ray $\rho_{56789}$ was not leaving $\operatorname{conv} v$ through the relative interior of $\operatorname{conv}\{v_1,v_{10}\}$,
we would get
\begin{equation}\label{eqeas3}
a_7<a_6
\end{equation}
as the analogue of~\eqref{eqeas2}.
Similarly, the analogue of~\eqref{eqeas1} for $\rho_{23458}$ reads
$$\frac{\dist(v_{8}, v_{4}\wedge v_{5})}{\dist(v_{3}, v_{4}\wedge v_{5})}<1.5$$
or $L_{56}\sin a_{5}+L_{67}\sin(a_5+a_6)+L_{78}\sin(a_5+a_6+a_7)<1.5L_{34}\sin a_{4}$, which implies
\begin{equation}\label{eqeas4}
(L_{56}+L_{78})\sin a_5<1.5L_{34}\sin a_{4}.
\end{equation}
Similarly, the assumption $\gamma_5\geqslant1$ in the formulation of the lemma means that
$$\frac{\dist(v_{8}, v_{6}\wedge v_{7})}{\dist(v_{5}, v_{6}\wedge v_{7})}\geqslant1$$
or $L_{78}\sin a_{7}\geqslant L_{56}\sin a_{6}$. Together with~\eqref{eqeas3}, this inequality implies
$L_{78}>L_{56}$, and we further apply~\eqref{eqeas4} and get $2L_{56}\sin a_{5}<1.5L_{34}\sin a_{4}$, or
$$
2\dist(v_{6}, v_{4}\wedge v_{5})<1.5\dist(v_{3}, v_{4}\wedge v_{5}),
$$
which means that $\gamma_3<3/4$ and contradicts to the assumption $\gamma_3\geqslant1$.
\end{proof}

\begin{defn}\label{defperf}
A correct sequence $(v_2,v_3,\ldots,v_k)$ is called \textit{perfect} if the quantity $\gamma_\tau$ as in~\eqref{eqgamma} satisfies $\gamma_\tau\geqslant1$ for all $\tau\in\{3,4,\ldots,k-3\}$.
\end{defn}

We proceed to the main result on the extension complexity of slanted sequences.

\begin{thr}\label{thrslant}
Let $\delta>1$ be an integer and $n=8\delta^2$. Let $v=(v_0,\ldots,v_n)$ be a $\pi/6$-thin sequence that is cw-slanted to the angle $\pi/3$ with tolerance $2\delta$. Then $v$ has a subsequence of size at least $0.25\delta^2$ and extension complexity at most $3\delta$.
\end{thr}

\begin{proof}
Assuming without loss of generality that $v$ is written in the cw-order, we build a cw-decreasing subsequence $u=(u_0, u_0',w_{01},u_1,u_1',w_{12},\ldots,u_h,u_h')$ of $v$. We take $i_0=0$, $u_0=v_0$, $u_0'=v_{4\delta}$, $w_{01}=v_{8\delta}$, and for any $q>0$, we inductively define
$$u_q=v_{4\delta i_q},\,\,\,\,u_q'=v_{4\delta(i_q+1)},\,\,\,\,w_{q\,q+1}=v_{4\delta(i_q+2)},$$
where $i_q\in\{i_{q-1}+3, i_{q-1}+4,\ldots,2\delta-4\}$ is the smallest index such that
$$\dist(v_{4\delta(i_q+2)},v_{4\delta i_q}\wedge v_{4\delta(i_q+1)})<\dist(w_{q-1\,q},v_{4\delta i_q}\wedge v_{4\delta(i_q+1)}),$$
and if no such $i_q$ exists, then we take $u_q=v_{n-4\delta}$, $u_q'=v_n$, and we set $h=q$. The situation splits into the two cases, which we need to treat separately.

\textit{Case 1.} Assume $h\geqslant 0.125\delta$. According to our definitions, the sequence $u$ is cw-decreasing, so we can apply Corollary~\ref{corcwdecr} with the subsequence $$(u_0,u_0',w_{01},u_1,u_1',w_{12},\ldots,u_\tau,u_\tau')$$ in the role of $v'$, where we can take $\tau=\lceil\delta/8\rceil$ because of the initial assumption of Case~1. This gives a subsequence of $v$ with at least $2\delta\tau\geqslant0.25\delta^2$ vertices and extension complexity not exceeding $2\delta+2\tau<3\delta$. 

\textit{Case 2.} Assume $h<0.125\delta$. For any $q\in\{0,\ldots,h-1\}$, we define the sequence
$$\sigma_q:=\left(v_{4\delta(i_q+2)},v_{4\delta(i_q+3)},\ldots,v_{4\delta(i_{q+1}+1)}\right),$$
which consists, in other words, of those points of the form $v_{4\delta\mathbb{Z}}$ that lie between $w_{q \,q+1}$ and $u_{q+1}'$. According to the choice of $i_q$ above, the sequence $\sigma_q$ is perfect in the sense of Definition~\ref{defperf}. For trivial reasons, every such $\sigma_q$ contains
$$m_q:=\left\lfloor\frac{|\sigma_q|}{8}\right\rfloor$$
disjoint subsequences of \textit{eight} consecutive vertices each. Any such subsequence is to be called \textit{main}, and every main subsequence is perfect because every $\sigma_q$ is perfect. The total number of main subsequences is $m:=m_0+\ldots+m_{h-1}$, and since they are disjoint, we have $m\leqslant(2\delta-1)/8<0.25\delta$. We also have
$$m=
\sum\limits_{q=0}^{h-1}\left\lfloor\frac{|\sigma_q|}{8}\right\rfloor\geqslant \sum\limits_{q=0}^{h-1}\left(\frac{|\sigma_q|}{8}-\frac{7}{8}\right)=\frac{1}{8}\left(\sum\limits_{q=0}^{h-1}|\sigma_q|\right)-\frac{7h}{8}=\frac{2\delta-1}{8}-\frac{7h}{8},$$
which implies $m\geqslant 0.25\delta-h$.

According to Lemma~\ref{lemout}, every main subsequence contains a further subsequence $\mu=(\mu_1,\mu_2,\mu_3,\mu_4,\mu_5)$ such that the ray $\rho(\mu,1,2,3,4,5)$ leaves $\conv v$ through the interior of $\conv\{v_0, v_n\}$. According to Definition~\ref{defslant} and Observation~\ref{obsslant1}, the ray $\rho(\mu',1,2,3',4,5)$ leaves $\conv v$ through the interior of $\conv\{v_0, v_n\}$ not only for $\mu_{3'}=\mu_3$ but also for at least $2\delta-1$ additional choices of a point $\mu_{3'}$ between $\mu_2$ and $\mu_3$ on the sequence $v$, where $\mu'=(\mu_1,\mu_2,\mu_{3'},\mu_4,\mu_5)$ denotes the sequence obtained from $\mu$ by replacing $\mu_3$ with a such vertex $\mu_{3'}$. According to Lemma~\ref{lemsameedge}, we obtain a $G$-good sequence if we take the concatenation of all $(\mu_1,\mu_2,\mu_{3'},\mu_4,\mu_5)$ in which $\mu_{3'}$ is arbitrarily chosen from the $2\delta$ vertices mentioned above (with $G$ being the set of the positions of the $\mu_{3'}$'s). Finally, we set $\mathcal{U}$ to be the sequence whose vertex set is the union of $\{\mu_1,\mu_2,\mu_4,\mu_5\}$ together with all possible choices of $\mu_{3'}$ over all sequences $\mu$ as above. We see that the sequence $\mathcal{U}$ has $$4m+2\delta m\geqslant (\delta-4h)+0.5\delta^2-2\delta h\geqslant (\delta-4\cdot 0.125\delta)+0.5\delta^2-2\delta\cdot0.125\delta\geqslant 0.25\delta^2$$
vertices, and, according to Lemma~\ref{lemgoodgood}, the extension complexity of $\mathcal{U}$ does not exceed $5m-m+2\delta\leqslant 3\delta$.
\end{proof}

\section{Completing the proof}\label{seccompl}

In this section, we put our technical results together and complete the proof of Theorem~\ref{thrmain}. First, we prove that a $\pi/6$-thin sequence admits a sufficiently large subsequence with small extension complexity.

\begin{thr}\label{thrthin}
Let $v$ be a sequence with $n=1024\tau^3+8\tau$ vertices, where $\tau$ is a positive integer. If $v$ is $\pi/6$-thin, then $v$ contains a subsequence with at least $4\tau^2$ vertices and extension complexity at most $12\tau$. 
\end{thr}

\begin{proof}
We apply Lemma~\ref{lemextrrob}, assuming $\beta=\pi/3$, $\delta=4\tau$, and $m=8\delta^2+1$.

\textit{Case 1.} The assertion of Lemma~\ref{lemextrrob} holds. Then $v$ contains a subsequence with $(6+\delta)\tau=4\tau^2+6\tau$ vertices and extension complexity not exceeding $6\tau+\delta+1=10\tau+1$, which satisfies the bounds desired in the current theorem.

\textit{Case 2.} The assumption of Lemma~\ref{lemextrrob} does not hold, which means that $v$ admits a subsequence $u=(u_1,\ldots,u_m)$ which is either cw-slanted or ccw-slanted to the angle $\pi/3$ with tolerance $2\delta$. According to Theorem~\ref{thrslant}, such a sequence $u$ has a subsequence of at least $0.25\delta^2=4\tau^2$ points and extension complexity at most $3\delta=12\tau$. These bounds match the desired conclusion as well.
\end{proof}

\begin{cor}\label{corthin}
Let $v$ be a sequence with $n>263\,000$ vertices. If $v$ is $\pi/6$-thin, then it contains a subsequence $u$ with at least $$\sqrt[3]{n^2}/36$$ vertices with extension complexity not exceeding $\sqrt[3]{72n/43}$.
\end{cor}

\begin{proof}
We take $$\tau=\left\lfloor\sqrt[3]{\frac{n}{1032}}\right\rfloor$$
and apply Theorem~\ref{thrthin} to arbitrary subsequence of $v$ containing $1024\tau^3+8\tau$ points. So we are able to extract a subsequence of $v$ with extension complexity at most
$$12\tau=12\left\lfloor\sqrt[3]{\frac{n}{1032}}\right\rfloor\leqslant12\sqrt[3]{\frac{n}{1032}}=\sqrt[3]{\frac{72n}{43}}$$
and containing at least $$4\tau^2=4\left(\left\lfloor\sqrt[3]{\frac{n}{1032}}\right\rfloor\right)^2>4\left(\sqrt[3]{\frac{n}{1032}}-1\right)^2>\frac{\sqrt[3]{n^2}}{36}$$ vertices, where the last inequality holds in the desired range of $n>263\,000$.
\end{proof}

\begin{cor}\label{corthinmain}
Let $v$ be a correct sequence with $n$ vertices. If $v$ is $\pi/6$-thin, then its extension complexity does not exceed $324\sqrt[3]{n^2/129}$.
\end{cor}

\begin{proof}
For $n\leqslant263\,000$, the result is trivial because the desired bound is greater than $n$. For $n>263\,000$, we proceed by induction and extract a subsequence $u$ as in Corollary~\ref{corthin}. So we see that the sequence $v$ splits into the disjoint union of two subsequences,

(1) one of which is $u$, and it has extension complexity not exceeding $\sqrt[3]{72n/43}$;

(2) the other one, $w$, consists of all those vertices in $v$ that are not in $u$.

\noindent The bound in Corollary~\ref{corthin} guarantees that $w$ has at most
$n-\sqrt[3]{n^2}/36$
vertices; by the inductive hypotesis we have
$$\xc(\conv w)\leqslant \frac{324}{\sqrt[3]{129}}\sqrt[3]{\left(n-\frac{\sqrt[3]{n^2}}{36}\right)^2},$$
which gives together with Lemma~\ref{lemeasy2}
$$\xc(\conv v)\leqslant\frac{324}{\sqrt[3]{129}}\sqrt[3]{\left(n-\frac{\sqrt[3]{n^2}}{36}\right)^2}+\sqrt[3]{\frac{72n}{43}}<\frac{324n}{\sqrt[3]{129n}},$$
where the last inequality is valid for any positive integer $n$.
\end{proof}

Now we can prove Theorem~\ref{thrmain}. According to Observation~\ref{lemextrthin}, the vertex set of any convex $n$-gon $P$ can be respresented as a disjoint union of twelve subsets $V_1,\ldots, V_{12}$ such that any subset with $|V_i|>2$ can be written as a $\pi/6$-thin sequence. We assume without loss of generality that each of the first $k$ subsets $V_1,\ldots, V_k$ have cardinality greater than $2$, and each of the remaining subsets has cardinality at most $2$. According to Lemma~\ref{lemeasy2}, we have
$$\xc(P)\leqslant\xc(\conv V_1)+\ldots+\xc(\conv V_k)+|V_{k+1}|+\ldots+|V_{12}|,$$
and, denoting $|V_i|$ by $n_i$, we apply Corollary~\ref{corthinmain} and obtain
\begin{equation}\label{eqrem}
\xc(P)\leqslant\frac{324}{\sqrt[3]{129}}\cdot\left(\sum\limits_{i=1}^{12}\sqrt[3]{n_i^2}\right)
\end{equation}
with $n_1+\ldots+n_{12}=n$. Since $t\to\sqrt[3]{t^2}$ is a concave function on $t>0$, the right-hand side of~\eqref{eqrem} attains its maximum if $n_1=\ldots=n_{12}=n/12$. This gives
$$\xc(P)\leqslant\frac{324}{\sqrt[3]{129}}\cdot\left(12\sqrt[3]{\frac{n^2}{12^2}}\right)<147\,\sqrt[3]{n^2},$$
which is the desired bound.

\section{Concluding remarks}\label{secrem}

We proved that $\wcc(n)\leqslant147\,n^{2/3}$, that is, every convex $n$-gon can be obtained as a linear projection of a higher-dimensional polytope with at most $147\,n^{2/3}$ facets. As outlined in Section~\ref{secresold}, this refutes expectations expressed in recent literature that $\wcc(n)$ should be close to a linear function. In view of our result, a weaker conjecture on lower bounds of polygon complexity could be proposed. 

\begin{conj}\label{conpc}
One has $\wcc(n)=\sqrt{n}\cdot\alpha(n)$ with unbounded $\alpha(n)$.
\end{conj}

Besides from being an adaptation of previous expectations, this conjecture has some further supporting evidence. We recall a result of Padrol~\cite{Pad} stating that 
\begin{equation}\label{EqPadrol}
\operatorname{wcc}(d,n) \geqslant 2\sqrt{dn-d}-d+1,
\end{equation}
where $\operatorname{wcc}(d,n)$ is the largest possible extension complexity of a polytope with $n$ vertices in a $d$-dimensional space. If Conjecture~\ref{conpc} were false, we would have $\wcc(n)\in O(\sqrt{n})$, which would mean that the bound~\eqref{EqPadrol} becomes asymptotically optimal when restricted to the case $d=2$. This is not quite expected, because even for slowly growing $d$ the function $\operatorname{wcc}(d,n)$ can grow much faster than $\sqrt{dn}$. In fact, a known result on the so-called \textit{correlation polytope}~\cite{KW} gives
\begin{equation}\label{eqdisc2}\operatorname{wcc}\left(m^2,2^m\right)\geqslant 1.5^m,\end{equation}
and the authors of~\cite{VGGT} make a further conjecture that implies $\operatorname{wcc}\left(m^2,2^m\right)=2^m$ for any positive integer $m$. The application of~\eqref{EqPadrol} gives a bound similar to~\eqref{eqdisc2} but with $1.5$ replaced by $\sqrt{2}$, which is much weaker in the asymptotical sense. As a further piece of evidence towards Conjecture~\ref{conpc}, we note that the method presented in this paper does not seem to allow an $O(\sqrt{n})$ upper bound for the worst-case $n$-gon complexity, but we were able to construct polygons whose optimal extended formulations require the construction as in Theorem~\ref{lemacutemult} already with $n=9$. However, as it is the case with many questions on lower bounds in combinatorial optimization, a proof of Conjecture~\ref{conpc} remains elusive.

\section{Acknowledgments}
My interest to Question~\ref{questwcc} developed after a fruitful workshop on \textit{Communication complexity, Linear optimization, and Lower bounds for the nonnegative rank of matrices}, which took place at \textit{Schloss Dagstuhl} in February, 2013. I would like to thank Nicolas Gillis and Fran\c{c}ois Glineur for personal discussions on the topic and all other participants for the information that I learned on Question~\ref{questwcc} and similar issues from their talks and open problem sessions. The first version of this paper~\cite{myfirstver}, which contained a proof that $\wcc(n)$ is $o(n)$, came in 2014 as a result of the effort I spent after this workshop, but the paper was not published. I would like to thank the organizers and participants of another workshop on \textit{Limitations of Convex Programming: Lower Bounds on Extended Formulations and Factorization Ranks} hosted by \textit{Schloss Dagstuhl} in 2015, and, in particular, I am grateful to Arnau Padrol Sureda for pointing me to~\cite{PP} and a further discussion on the topic of Question~\ref{questwcc}.  Also, I would like to thank Nicolas Gillis and Arnaud Vandaele for a further interesting discussion that held during my visit to the \textit{University of Mons} in 2017, and I am grateful to Nicolas for invitation.

\end{document}